\documentclass[twoside]{amsart}

\usepackage{array}
\usepackage{color}             
\usepackage{amsfonts}
\usepackage{amssymb, amsmath}   
\usepackage{bm}
\usepackage[mathscr]{eucal}
\usepackage{tikz-cd}            
\usepackage{multicol}

\allowdisplaybreaks

\def\qb{\hfill $\Box$}

\def\cf{{\it cf.}\ }

\begin{document}

\title{On the $\mathfrak{S}_{q} \overleftrightarrow{\times} \mathfrak{S}_q$-invariant rational cohomology of the pure braid group $P_{2q}$}
\author{Yan Fu}
\author{Gefei Wang}

\date{}
\thanks{The first named author is supported by NSFC grant No.12271183. The second named author is partially supported by National Key R\&D Program of China 2020YFA0712800 and NSFC grant No.12131009.}
\address{Yan Fu, Chern Institute of Mathematics,
Nankai University, Tianjin, P. R. China}
\email{1120230012@mail.nankai.edu.cn}
\address{Gefei Wang, School of Mathematical Science,
	Peking University, Beijing,
	P. R. China}
\email{wanggefei@math.pku.edu.cn}

\begin{abstract} Let $\mathfrak{S}_{q} \overleftrightarrow{\times} \mathfrak{S}_q$ be the $\mathbb{Z}/2$-extension of the product of two symmetric groups $\mathfrak{S}_{q} \times \mathfrak{S}_q$. In this paper, we compute the $\mathfrak{S}_{q} \overleftrightarrow{\times} \mathfrak{S}_q$-invariant part of the  rational cohomology of the pure braid group $P_{2q}$.   This includes the rational cohomology of a spin hyperelliptic mapping class group of genus $g$, where $g+1=q$.
\end{abstract}

\maketitle

\section{Introduction}

Let $\Sigma_g$ be a compact connected oriented surface of genus 
$g$.
A spin structure $c$ on $\Sigma_g$
is a quadratic refinement of the intersection form on $H_1(\Sigma_g , \mathbb{Z}/2)$, characterized by the identity
$$c(x+y)=c(x)+c(y)+x \cdot y, \forall x,y \in H_1(\Sigma_g , \mathbb{Z}/2)$$
where $x \cdot y$ denotes the algebraic intersection number. The mapping class group acts on spin structures via its action on homology, and the spin mapping class group  associated with a spin structure $c$ is defined as  the isotropy group at  $c$ ({\it cf.} \cite{J80}).

The bulk of existing results on the  homology and cohomology of spin mapping class groups focus on the stable range ({\it cf.} e.g., \cite{RM12}, \cite{RM14}). However, there are unstable aspects in the cohomology of the spin mapping class group.  Our work focuses on revealing the unstable aspect in the cohomology of spin hyperelliptic mapping class groups.

Let $\mathcal{S}(\Sigma_g)$ be the hyperelliptic mapping class group for a compact connected oriented surface $\Sigma_g$ of genus $g$.   
The spin hyperelliptic mapping class group $\mathcal{S}(\Sigma_g ; c)$ is the isotropy group of $\mathcal{S}(\Sigma_g)$ at each spin structure c. The action of the group $\mathcal{S}(\Sigma_g)$ on the set of spin structures on the surface $\Sigma_g$ has $\lceil \frac{g}{2}\rceil+1$ orbits. Hence we have  $\lceil \frac{g}{2}\rceil+1$ spin hyperelliptic mapping class groups for $\Sigma_g$. The isotropy group $\mathfrak{G}_i$ of each orbit can be described by the direct product of two symmetric groups $\mathfrak{S}_{g+1+2i} \times \mathfrak{S}_{g+1-2i}$ or a $\mathbb{Z}/2$-extension of the direct product of two symmetric groups $\mathfrak{S}_{g+1} \overleftrightarrow{\times} \mathfrak{S}_{g+1}$ ({\it cf.} e.g., \cite{M84} Proposition 3.1, p.3.95, or \cite{W23}). 

Now we consider the rational cohomology groups $H^*(\mathcal{S}(\Sigma_{g};c))$.  From \cite{BH71} and since the Lyndon-Hochschild-Serre spectral sequence $E_2$-page equals $E_\infty$-page, \\
 $H^*(\mathcal{S}(\Sigma_{g};c))$ is isomorphic to the  $\mathfrak{G}_i$-invariant part of the rational cohomology group of the pure sphere braid group $SP_{2g+2}$, denoted by $H^*(SP_{2g+2})^{\mathfrak{G}_i}$. Moreover, for any $n \geq 2$, the natural surjection of the pure Artin braid group $P_n$ onto the pure sphere braid group $SP_n$ induces an injection of the rational cohomology groups (\cf e.g., \cite{K90}). Hence
the  $\mathfrak{G}_i$-invariant part of the rational cohomology group of  $P_{2g+2}$, denoted by $H^*(P_{2g+2})^{\mathfrak{G}_i}$, includes $H^*(\mathcal{S}(\Sigma_{g};c))$ as a subalgebra. 

Our principal objective is to establish a comprehensive computation of  $H^*(P_{n})^{\mathfrak{G}}$ for every $\mathfrak{G}$ such that  $H^*(\mathcal{S}(\Sigma_{g};c))$ is included as a subalgebra when $n=2g+2$.  Building upon the result of Lehrer and Solomon In \cite{LS86}, the second named author constructed a basis of the underlying graded $\mathbb{Q}$-vector space $H^{*}(P_n)^{\mathfrak{G}}$ and calculated $H^{*}(P_n)^{\mathfrak{G}}$  for the cases  $\mathfrak{G}=\mathfrak{S}_{n-q} \times \mathfrak{S}_{q}$, $n-q\geq q$ in \cite{W24}.
Crucially, the case
$\mathfrak{G}=\mathfrak{S}_{q} \overleftrightarrow{\times} \mathfrak{S}_q$, $n=2q$, remained unresolved. 

In this paper, we resolve the previously outstanding case by calculating $H^{*}(P_n)^{\mathfrak{G}}$ with  $\mathfrak{G}=\mathfrak{S}_{q} \overleftrightarrow{\times} \mathfrak{S}_q$, $n=2q$.
In section 2, we give some revisions on the  $(\mathfrak{S}_{q} \times \mathfrak{S}_{q})$-invariant part of the rational cohomology of  $P_{n}$, where $2q=n$.
In section 3, based on the fact that  $\mathfrak{S}_{q} \times \mathfrak{S}_q$ is a normal subgroup of
$\mathfrak{G}$, we consider a surjection 
$$\phi: H^{*}(P_n)^{\mathfrak{S}_{q} \times \mathfrak{S}_q}\rightarrow H^{*}(P_n)^{\mathfrak{G}}$$ of $\mathbb{Q}\mathfrak{G}$-modules and give the kernel of $\phi$.
 Since $H^{*}(P_n)^{\mathfrak{S}_{q} \times \mathfrak{S}_q}$ has already been calculated, we get  $H^{*}(P_n)^{\mathfrak{G}}$. In section 4, we  prove that

{\noindent\bf Theorem 4.2.} {\it $$\text{dim}(H^{*}(P_n)^{\mathfrak{G}})=\frac{\text{dim}(H^{*}(P_n)^{\mathfrak{S}_{q} \times \mathfrak{S}_q})+|E^P|}{2}-|K^P|.$$}
Where $E^P$ and $K^P$ are subsets of the basis of $H^{*}(P_n)^{\mathfrak{S}_{q} \times \mathfrak{S}_q}$.
We also derive a formula to compute the cardinality $|E^P|$ and $|K^P|$.

{\noindent\bf Acknowledgments} 

The second named author would like to thank Nariya Kawazumi for helpful discussions. All authors contribute equally.

\section{Revisions on $H^{*}(P_n)^{\mathfrak{S}_{q} \times \mathfrak{S}_q}$}
 By definition,  the pure braid group $P_n$ is the fundamental group of the configuration space
 \[X_n := \{(u_1, \cdots, u_n) \in \mathbb{C}^n\mid u_i \neq u_j \text{ for } i\neq j\}.\]
 From the study by Arnold \cite{A69}, the cohomology ring $H^*(P_{n})$ is the quotient of an exterior graded ring, denoted by $A(n)$, generated by
 $\left(
 \begin{array}{c}
 	n \\
 	2 \\
 \end{array}
 \right)$ differential forms $\omega_{i,j}=\omega_{j,i}=\frac{du_i-du_j}{u_i-u_j}$, $u_i \in \mathbb{C}$, $1\leq i <j \leq n$, modulo the $\left(
 \begin{array}{c}
 	n \\
 	3 \\
 \end{array}
 \right)$ relations
 \[\omega_{k.l}\omega_{l,m}+\omega_{l,m}\omega_{m,k}+\omega_{m,k}\omega_{k,l}=0,\]
 where $1\leq k< l < m\leq n$.
 
 The symmetric group $\mathfrak{S}_n$ acts on $A(n)$ in a natural way.
 
 For a partition $\lambda=(\lambda_{1}\geq\lambda_{2}\geq\cdots\geq\lambda_{j}>0)$ of $n=\lambda_{1}+\lambda_{2}+\cdots+\lambda_{j}$, let $c_{\lambda}=c_{\lambda_1}c_{\lambda_2}\cdots c_{\lambda_j}$ be a element of the symmetric group $\mathfrak{S}_n$, where
 $$c_{\lambda_i}=(\lambda_1+\lambda_2+\cdots+\lambda_{i-1}+1,\lambda_1+\lambda_2+\cdots+\lambda_{i-1}+2,\cdots,\lambda_1+\lambda_2+\cdots+\lambda_{i})$$ is a cycle of length $\lambda_i$.
 Let $N_{\lambda}$ be the subgroup of $\mathfrak{S}_n$ generated by $\{\nu_i|\lambda_i=\lambda_{i+1}\}$, where
 \begin{align*}
 	\nu_i=&(\lambda_1+\lambda_2+\cdots+\lambda_{i-1}+1,\lambda_1+\lambda_2+\cdots+\lambda_{i}+1)\\
 	&(\lambda_1+\lambda_2+\cdots+\lambda_{i-1}+2,\lambda_1+\lambda_2+\cdots+\lambda_{i}+2)\\
 	&\cdots\\
 	&(\lambda_1+\lambda_2+\cdots+\lambda_{i},\lambda_1+\lambda_2+\cdots+\lambda_{i+1}),
 \end{align*}
 $C_{\lambda}=\langle c_{\lambda_1}\rangle\times \langle c_{\lambda_2}\rangle\times\cdots\times \langle c_{\lambda_j}\rangle$  the direct product of cyclic groups $\langle c_{\lambda_i}\rangle$.
 Let $Z_{\lambda}$ be the centralizer of the element $c_{\lambda}$.
 Then the  centralizer $Z_{\lambda}$ of the element $c_{\lambda}$ is the semi-product of $N_{\lambda}$ and $C_{\lambda}$, $Z_{\lambda}=N_{\lambda}\ltimes C_{\lambda}$.
 
 Assume that $\alpha_{\lambda}$ is the character of $N_{\lambda}$ defined by $\alpha_{\lambda}(\nu_i)=(-1)^{\lambda_i+1}$, $\phi_{\lambda}=(\phi_{\lambda_1}\otimes\phi_{\lambda_2}\otimes\cdots\otimes\phi_{\lambda_{j}})\varepsilon$ is the character of $C_{\lambda}$, where $\phi_{\lambda_i}(g_i)=e^{\frac{2\pi i}{\lambda_i}}$ and $\varepsilon$ is the sign character of $\mathfrak{S}_n$.
 From \cite{LS86}, the $i$-th cohomology group of $P_n$ has a $\mathbb{C}\mathfrak{S}_n$-module isomorphism
 \[ H^i(P_{n},\mathbb{C}) \cong \bigoplus_{\lambda}\text{Ind}_{Z_{\lambda}}^{\mathfrak{S}_n}(\zeta_{\lambda}),\]
 where $\lambda=(\lambda_{1}\geq\lambda_{2}\geq\cdots\geq\lambda_{j}>0)$ are those $j$ 
 partitions of $n$ such that $i+j=n$ and $\zeta_{\lambda}=\alpha_{\lambda}\phi_{\lambda}$ is the 1-dimensional linear representation of $Z_\lambda$.
 
  Let $\mathfrak{G}$ be a subgroup of $\mathfrak{S}_{n}$.
 According to the classical linear representation theory, there exists a $\mathbb{C}\mathfrak{G}$-module isomorphism \[ H^*(P_{n},\mathbb{C}) \cong \bigoplus_{\lambda}\text{Res}_{\mathfrak{G}}\text{Ind}_{Z_{\lambda}}^{\mathfrak{S}_n}(\zeta_{\lambda})
 \cong\bigoplus_{\lambda}\bigoplus_{s\in \mathfrak{G}\backslash\mathfrak{S}_n / Z_{\lambda}}\text{Ind}_{(Z_{\lambda})_s}^{\mathfrak{G}}(s \otimes \zeta_{\lambda}),\]
 where $(Z_{\lambda})_s=sZ_{\lambda}s^{-1} \cap \mathfrak{G}$ and $s \otimes \zeta_{\lambda}$ is defined by $(s \otimes \zeta_{\lambda})(g)=\zeta_{\lambda}(s^{-1}gs)$.
 Hence one can compute the $\mathfrak{G}$-invariant part of $H^*(P_{n})$ by 
 $$(s \otimes \zeta_{\lambda},1)_{(Z_{\lambda})_s}=\frac{1}{|(Z_{\lambda})_s|}\sum_{g \in(Z_{\lambda})_s}(s \otimes \zeta_{\lambda})(g)\overline{1}(g)=\frac{1}{|(Z_{\lambda})_s|}\sum_{g \in(Z_{\lambda})_s}(s \otimes \zeta_{\lambda})(g),$$
 where $\overline{1}=1$ is the $1$-dimensional trivial representation.
 More explicitly, by Corollary $2.2$ in \cite{W24}, there is  an isomorphism of $\mathbb{Q}$-vector spaces
 $$H^*(P_{n})^{\mathfrak{G}}\cong\bigoplus_{\lambda}\bigoplus_{(s \otimes \zeta_{\lambda},1)_{(Z_{\lambda})_s}=1 \atop s\in \mathfrak{G}\backslash\mathfrak{S}_n / Z_{\lambda}}(s \otimes \zeta_{\lambda}).$$

  Now  let $n=2q$ and $\mathfrak{G}=\mathfrak{S}_{q} \times \mathfrak{S}_q$, $s\in \mathfrak{G}\backslash\mathfrak{S}_{n} / Z_{\lambda}$. To calculate all $s\in \mathfrak{G}\backslash\mathfrak{S}_n / Z_{\lambda}$ such that $(s \otimes \zeta_{\lambda},1)_{(Z_{\lambda})_s}=1$, it is necessary to establish a complete invariant of $\mathfrak{G}\backslash\mathfrak{S}_n / Z_{\lambda}$.
  
{\noindent\bf Definition 2.1.}
{\it Let $(\mathbb{Z}/2)^{n}$ be the set of functions $\delta:[2q]\rightarrow\mathbb{Z}/2$ and $s\in(\mathfrak{S}_{q}\times\mathfrak{S}_q)\backslash\mathfrak{S}_{n}$. We define $\delta_s\in(\mathbb{Z}/2)^{n}$ by $$\delta_s(a)=\left\{\begin{array}{rl}0 ,& \text{if }1\leq s(a)\leq q,\\
		1 , & \text{if }q+1\leq s(a) \leq n.\end{array}\right.$$}

Recall that $\lambda=(\lambda_1\geq\lambda_2\geq\cdots\geq \lambda_j>0)$ is a partition of $n$ of length $j$.

{\noindent\bf Definition 2.2.}(\cite{W24} Definition 2.7.)
{ \it Let
	\begin{align*}B_{i}(\delta_s)=&\{\lambda_1+\cdots+\lambda_{i-1}+1\leq a\leq\lambda_{1}+\cdots+\lambda_{i}|\delta_s(a)=1\} \\
		\subset &\{\lambda_{1}+\cdots+\lambda_{i-1}+1,\cdots,\lambda_{1}+\cdots+\lambda_{i}\}.
	\end{align*}
	We  write  $B_{i}(\delta_s)=\{b_1<\cdots<b_{d_i}\}$ with $d_i=|B_{i}(\delta_s)|$.  We denote $\langle c_{d_i}\rangle$ the cyclic group generated by the cycle $c_{d_i}$ of length $d_i$.
	If $|B_{i}(\delta_s)|\geq 2$, we define an element in $(\mathbb{Z}_{\geq0})^{d_i}/\langle c_{d_i}\rangle$, denoted by $\chi_i(\delta_s)$, by
	$$
	\chi_i(\delta_s)=(b_2-b_1-1,\cdots, b_{d_i}-b_{d_i-1}-1, \lambda_i-b_{d_i}+b_{1}-1),
	$$ 
	while we define $\chi_i(\delta_s)=(\lambda_i-1)$ if $|B_{i}(\delta_s)|=1$ and $\chi_i(\delta_s)=\emptyset_{\lambda_i}$ if $B_{i}(\delta_s)=\emptyset$. We call $\chi_i(\delta_s)$ the $i$-th invariant cycle of $\delta$. We define a  multiset $\chi(\delta_s)$ 
	by
	$$\chi(\delta_s)=\{\chi_1(\delta_s),\chi_2(\delta_s),\cdots,\chi_{j}(\delta_s)\},$$ where
	$$\chi_i(\delta_s)\in \left(\{\emptyset_{\lambda_i} \}\cup(\bigcup_{d\geq 1}(\mathbb{Z}_{\geq0})^{d}/\langle c_{d}\rangle)\right),$$
	$1\leq i \leq j$ and call it the full invariant cycle set of $\delta_s$.
	We define a lexicographic order in the multiset $\chi(\delta_s)$ by
	\begin{enumerate}
		\item We choose a canonical representative of 
		$\chi_{i}(\delta_s)$
		that is lexicographically minimal.
		
		\item $\chi_i(\delta_s)\geq\emptyset_{\lambda_i}$.
		
		\item If $\lambda_i>\lambda_{i+1}$, then $\chi_i(\delta_s)>\chi_{i+1}(\delta_s)$.
		
		\item  If $\lambda_i=\lambda_{i+1}$, then $|B_i(\delta_s)|\geq |B_{i+1}(\delta_s)|$.
		
		\item If $\lambda_i=\lambda_{i+1}$ and $|B_i(\delta_s)|= |B_{i+1}(\delta_s)|$, then  $\chi_i(\delta_s)\leq\chi_{i+1}(\delta_s)$, where the order is given by $(1)$.
	\end{enumerate}
}

 From \cite{W24} Section 2, $\chi(\delta_s)$ is a complete invariant of $s\in\mathfrak{G}\backslash \mathfrak{S}_{n}\slash Z_{\lambda}$.

 Consider the full invariant cycle set 
$$\chi(\delta_{s})=\{\chi_1(\delta_{s}),\cdots,\chi_j(\delta_{s})\}$$
of an element $s\in\mathfrak{G}\backslash \mathfrak{S}_{n}\slash Z_{\lambda}$.
We calculate those $s\in\mathfrak{G}\backslash \mathfrak{S}_{n}\slash Z_{\lambda}$ such that  $(s \otimes \zeta_{\lambda},1)_{(Z_{\lambda})_s}=1$ by using $\chi(\delta_s)$. We denote $\langle c_{d_i}\rangle_{\chi_i(\delta_s)}$ the isotropy group of the $\langle c_{d_i}\rangle$ action on  $\chi_{i}(\delta_{s})$. 	

{\noindent\bf Theorem 2.3.}(\cite{W24} Theorem 2.8.)
{\it $(s \otimes \zeta_{\lambda},1)_{(Z_{\lambda})_s}=1$ if and only if $\chi(\delta_s)$ satisfies
	\begin{enumerate}
		\item $1\leq |B_i(\delta_s)|\leq \lambda_i-1$ for every $\lambda_i\geq 3$.
		\item If $\lambda_i\geq 3$ with $\lambda_i\equiv0,1,3$ mod $4$, then $\langle c_{d_i}\rangle_{\chi_i(\delta_s)}$ is of order $1$.
		\item If $\lambda_i\geq 2$ with $\lambda_i\equiv2$ mod $4$, then $\langle c_{d_i}\rangle_{\chi_i(\delta_s)}$ is at most of order $2$.
		\item If there exist some $\lambda_i=\lambda_{i+1}=2k$,  then $|B_i(\delta_{s})|>|B_{i+1}(\delta_{s})|$ or $|B_i(\delta_{s})|= |B_{i+1}(\delta_{s})|$ with $\chi_i(\delta_{s})<\chi_{i+1}(\delta_{s})$.
\end{enumerate}}
 	
  We denote by $\Lambda(n-*)$ the set of all  partitions $\lambda$ of $n$ of length $(n-*)$.

{\noindent\bf Definition 2.4.}(\cite{W24} Definition 3.1.)
{\it Let $\lambda=(\lambda_1,\lambda_2,\cdots,\lambda_a,1,1,\cdots)\in\Lambda(n-*)$,    $\lambda_i\geq 2$, $1\leq i \leq a$. We define $$d=(d_1,d_2,\cdots,d_a)\in(\mathbb{Z}_{\geq0})^a$$ and
	$$(\lambda,d)=\{(\lambda_1,d_1),(\lambda_2,d_2),\cdots,(\lambda_a,d_a)\}$$ which satisfies $d_i\geq d_{i+1}$ if $\lambda_i=\lambda_{i+1}$. If $d$ satisfies
	$d_1+d_2+\cdots+d_a\leq q$, then it is easy to see that there exists some  $\chi_i(\delta_s)$ such that $|B_i(\delta_s)|=d_i$.
	We define
	$$\Pi(\lambda_i,d_i)=\{\chi_i(\delta_s)\in\left(\{\emptyset_{\lambda_i}\}\cup((\mathbb{Z}_{\geq0})^{d_i}/\langle c_{d_i}\rangle)\right)||B_i(\delta_s)|=d_i,(s\otimes\phi_{\lambda_i}\varepsilon,1)_{(Z_{\lambda})_s}=1\}$$
	and denote by $|\Pi(\lambda_i,d_i)|$ the cardinality of $\Pi(\lambda_i,d_i)$.
}  
 
 According to the definition of $\chi$ and from Theorem 2.3, if we take $$\bigcup_{\lambda}\{s\otimes\zeta_{\lambda}|(s\otimes\zeta_{\lambda},1)_{(Z_{\lambda})_s}=1\}$$ as a basis of $H^*(P_{n})^{\mathfrak{G}}$, then we can get each basis element $s\otimes\zeta_{\lambda}$ by $\chi_i(\delta_s)\in\Pi(\lambda_i,d_i)$ for $1\leq i \leq a$ and some $d_i$. ({\it cf.} \cite{W24} Theorem 2.9 and Theorem 3.2)

\section{the $\mathfrak{S}_{q} \overleftrightarrow{\times} \mathfrak{S}_q$-invariant part of $H^{*}(P_n)$}

In this section, we focus on the $\mathbb{Z}/2$-extension of the direct product of two symmetric groups $\mathfrak{G}=\mathfrak{S}_{q} \overleftrightarrow{\times} \mathfrak{S}_q$, where the $\mathbb{Z}/2$-extension is given by $$\sigma=(1,2q)(2,2q-1)\cdots(q,q+1).$$

For a fixed invariant cycle $\chi_i(\delta_s)$ defined in Definition $2.2$, we consider the invariant cycle $\chi_i(\delta_{\sigma s})$ at first.
 It is easy to see that if $\chi_i(\delta_s)=(0,0,\cdots,0)$, then $\chi_i(\delta_{\sigma s})=\emptyset_{\lambda_i}$.

{\noindent\bf Lemma 3.1.}
{\it Let  $$\chi_i(\delta_s)=(\underbrace{0,\cdots,0}_{k_1},a_1,\underbrace{0,\cdots,0}_{k_2},a_2,\cdots,\underbrace{0,\cdots,0}_{k_m},a_m)$$
 be an invariant cycle defined in Definition $2.2$, $k_i\geq0$, $a_i\geq1$, $1\leq i\leq m$. Then
$$\chi_i(\delta_{\sigma s})=(k_1+1,\underbrace{0,\cdots,0}_{a_1-1},k_2+1,\underbrace{0,\cdots,0}_{a_2-1},\cdots,k_m+1,\underbrace{0,\cdots,0}_{a_m-1}).$$
Furthermore, there exists a bijection between $\Pi(\lambda_i,d_i)$ and $\Pi(\lambda_i,\lambda_i-d_i)$.
}

{\noindent\bf Proof.} Let
$$\chi_i(\delta_s)=(\underbrace{0,\cdots,0}_{k_1},a_1,\underbrace{0,\cdots,0}_{k_2},a_2,\cdots,\underbrace{0,\cdots,0}_{k_m},a_m).$$
From Definition $2.2$, we consider
$B_i(\delta_s)=\{b_1<\cdots<b_{d_i}|\delta_s(b_j)=1, 1\leq j \leq d_i\}$ with $d_i=|B_i(\delta_s)|$.
Then we have 

\begin{align*}
	\lambda_i&=k_1+\cdots+k_m+m+a_1+\cdots+a_m,\\
  d_i&=k_1+\cdots+k_m+m,\\
  b_{1+j}&=b_1+j, 1\leq j \leq k_1,\\
  b_{2+k_1}&=b_{1+k_1}+a_1+1,\\
  b_{2+k_1+j}&=b_{2+k_1}+j, 1\leq j \leq k_2,\\
  b_{3+k_1+k_2}&=b_{2+k_1+k_2}+a_2+1,\\
  &\vdots\\
  b_{m+k_1+\cdots+k_{m-1}+j}&=b_{m+k_1+\cdots+k_{m-1}}+j, 1\leq j \leq k_m.
\end{align*}
 By the definition of $\sigma$, we have $$\delta_{\sigma s}(i)=\delta(\sigma(s(i)))=\delta(2q+1-s(i))=\left\{\begin{array}{rl}1 ,& \text{if }1\leq s(i)\leq q,\\
0 , & \text{if }q+1\leq s(i) \leq 2q.\end{array}\right.$$
Since $\delta_{s}(b_j)=1$, $1\leq j\leq d_i$, $b_j\in B_i(\delta_s)$, we have $\delta_{\sigma s}(b_j)=0$, $1\leq j\leq d_i$, $b_j\in B_i(\delta_s)$. Thus we have
\begin{align*}
	B_i(\delta_{\sigma s})&=\{\lambda_{1}+\cdots+\lambda_{i-1}+1,\cdots,\lambda_{1}+\cdots+\lambda_{i}\}\setminus B_i(\delta_{s}),\\
	&=\left\{\begin{array}{c}
		\lambda_{1}+\cdots+\lambda_{i-1}+1<\lambda_{1}+\cdots+\lambda_{i-1}+2<\cdots<b_1-1\\
	<b_{1+k_1}+1<b_{1+k_1}+2<\cdots<b_{1+k_1}+a_1\\
	<b_{2+k_1}+k_2+1<b_{2+k_1}+k_2+2<\cdots<b_{2+k_1}+k_2+a_2\\
	\vdots\\
	<b_{m+k_1+\cdots+k_{m}}+1<b_{m+k_1+\cdots+k_{m}}+2<\cdots<\lambda_{1}+\cdots+\lambda_{i}
	\end{array}\right\},\\
	|B_i(\delta_{\sigma s})|&=\lambda_i-d_i=a_1+\cdots+a_m.
\end{align*}
Hence we have $$\chi_i(\delta_{\sigma s})=(k_1+1,\underbrace{0,\cdots,0}_{a_1-1},k_2+1,\underbrace{0,\cdots,0}_{a_2-1},\cdots,k_m+1,\underbrace{0,\cdots,0}_{a_m-1}),$$
where we remark that
\begin{align*}
	k_1+1&=b_{1+k_1}+1-(b_1-1)-1,\\ k_j+1&=b_{j+k_1+\cdots+k_{j-1}}+k_j+1-(b_{j-1+k_1+\cdots+k_{j-1}}+a_{j-1})-1, 2\leq j \leq m.
\end{align*}

Now assume that $\chi_i(\delta_s)\in \Pi(\lambda_i,d_i)$. We have $|B_i(\delta_{\sigma s})|=\lambda_i-d_i$. Note that if $\langle c_{d_i}\rangle_{\chi_i(\delta_s)}$ is of order $p$, then $p\mid d_i$ and $p\mid m$. Consider $\chi_i(\delta_s)$,  we have $$\sum_{i=1}^{\frac{m}{p}}(k_i+1)=\frac{d_i}{p},$$ $k_i=k_{i+\frac{jm}{p}}$ and $a_i=a_{i+\frac{jm}{p}}$ for $1\leq i \leq \frac{m}{p}$, $1\leq j \leq p-1$. Thus for $\chi_i(\delta_{\sigma s})$, we get
$$\sum_{i=1}^{\frac{m}{p}}(a_i-1+1)=\frac{\lambda_i-d_i}{p}$$
and
$p\mid \lambda_i-d_i$.
These show that the order of $\langle c_{d_i}\rangle_{\chi_i(\delta_s)}$ equals the order of $\langle c_{\lambda_i-d_i}\rangle_{\chi_i(\delta_{\sigma s})}$. From Theorem $2.3$, we have $\chi_i(\delta_{\sigma s})\in \Pi(\lambda_i,\lambda_i-d_i)$.

These prove the Lemma.\qb

{\noindent\bf Definition 3.2.} {\it Let $\chi_i(\delta_s)=(a_1,a_2,\cdots,a_d)$ be an invariant cycle defined in Definition $2.2$.
We define an invariant cycle associated with  $\chi_i(\delta_s)$, denoted by $\overleftarrow{\chi_i(\delta_s)}$, by $$\overleftarrow{\chi_i(\delta_s)}=(a_d,\cdots,a_2,a_1).$$
We define the $\sigma$ action on $\chi_i(\delta_s)$  by 
 $$\sigma\circ\chi_i(\delta_s)=\overleftarrow{\chi_i(\delta_{\sigma s})}.$$
We denote $\overline{\chi_i}(\delta_s)=\sigma\circ\chi_i(\delta_{s})$ . Then $\sigma$ acts on the full invariant cycle set
$$\chi(\delta_s)=\{\chi_1(\delta_s),\chi_2(\delta_s),\cdots,\chi_j(\delta_s)\}$$ in a natural way by
$$\sigma\circ\chi(\delta_s)=\{\overline{\chi_1}(\delta_s),\overline{\chi_2}(\delta_s),\cdots,\overline{\chi_j}(\delta_s)\}.$$
We also denote $\overline{\chi}(\delta_s)=\sigma\circ\chi(\delta_s)$.}

Since $$\sigma^2\circ\chi_i(\delta_s)=\sigma\circ\overleftarrow{\chi_i(\delta_{\sigma s})}=\overleftarrow{\overleftarrow{\chi_i(\delta_{s})}}=\chi_i(\delta_{s}),$$ the $\sigma$ action on $\chi_i(\delta_{s})$ is well-defined.

Note that $\mathfrak{S}_{q} \times \mathfrak{S}_q$ is a normal subgroup of
$\mathfrak{G}$, $\sigma$ acts on $H^{*}(P_n)$ in a natural way. On the other hand,  $ H^{*}(P_n)^{\mathfrak{S}_{q} \times \mathfrak{S}_q}$ and  $H^{*}(P_n)^{\mathfrak{G}}$ can be determined by the averaging operators  $$\mathcal{A}_{\mathfrak{S}_{q} \times \mathfrak{S}_q}=\frac{1}{|\mathfrak{S}_{q} \times \mathfrak{S}_q|}\sum_{\tau\in\mathfrak{S}_{q} \times \mathfrak{S}_q}\tau$$ and
$$\mathcal{A}_{\mathfrak{G}}=\frac{1}{|\mathfrak{G}|}\sum_{\tau\in\mathfrak{G}}\tau$$
act on $H^{*}(P_n)$. 
Hence
there exists a commutative diagram of $\mathbb{Q}\mathfrak{G}$-modules
\[ 
\begin{tikzcd}
	H^{*}(P_n) \arrow[d, two heads, "\mathcal{A}_{\mathfrak{S}_{q} \times \mathfrak{S}_q}"]  \arrow[dr, two heads, "2\mathcal{A}_{\mathfrak{G}} "]	
	\\
	H^{*}(P_n)^{\mathfrak{S}_{q} \times \mathfrak{S}_q} \arrow[r, two heads, "1+\sigma"] & 
	H^{*}(P_n)^{\mathfrak{G}}
\end{tikzcd}
.\]
We denote by $\phi$ the surjection from $H^{*}(P_n)^{\mathfrak{S}_{q} \times \mathfrak{S}_q}$ to $	H^{*}(P_n)^{\mathfrak{G}}$.
More explicitly, assume that $s\otimes\zeta_{\lambda}\in H^{*}(P_n)^{\mathfrak{S}_{q} \times \mathfrak{S}_q}$, then
$$\phi(s\otimes\zeta_{\lambda})=s\otimes\zeta_{\lambda}+\sigma\circ (s\otimes\zeta_{\lambda}).$$ 

In \cite{W24} Theorem $2.9$, for a fixed partition $\lambda=(\lambda_1,\lambda_2,\cdots,\lambda_a,1,\cdots)$, the second named author gave a ring structure of $H^{*}(P_n)^{\mathfrak{S}_{q} \times \mathfrak{S}_q}$ by considering a surjection
\[   \bigotimes_{i=1}^{a}H^{\lambda_i-1}(P_n)^{\mathfrak{S}_{q} \times \mathfrak{S}_q}
\rightarrow H^{*}(P_n)^{\mathfrak{S}_{q} \times \mathfrak{S}_q}.
\]
 If we consider  a sequence $$(\chi_1(\delta_s),\chi_2(\delta_s),\cdots,\chi_a(\delta_s))$$
 such that $$\sum_{i=1}^a|B_i(\delta_s)|\leq q,$$
 then it is associated with an element in $\bigotimes_{i=1}^{a}H^{\lambda_i-1}(P_n)^{\mathfrak{S}_{q} \times \mathfrak{S}_q}$ such that $$\{\chi_i(\delta_s),(0),\cdots,(0),\emptyset_1,\cdots,\emptyset_1\}$$ is associated with an element in $H^{\lambda_i-1}(P_n)^{\mathfrak{S}_{q} \times \mathfrak{S}_q}$, $1\leq i\leq a$. By the surjection, up to sign considerations,  we get an element in $ H^{*}(P_n)^{\mathfrak{S}_{q} \times \mathfrak{S}_q}$ whose full invariant cycle set is $$\{\chi_1(\delta_s),\chi_2(\delta_s),\cdots,\chi_a(\delta_s),(0),\cdots,(0),\emptyset_1,\cdots,\emptyset_1 \},$$ where the sign is depending on the rearrangement of the sequence $$(\chi_1(\delta_s),\chi_2(\delta_s),\cdots,\chi_a(\delta_s))$$
 by the lexicographic order we defined in Definition 2.2.
 
 We remark that for a fixed $\lambda=\{\lambda_1\geq\lambda_2\geq\cdots\}\in\Lambda(n-*)$,
 $\chi_i(\delta_s)=\overline{\chi_i}(\delta_s)$  only if $\lambda_i=2d_i$. 
 The $\langle c_{d_i}\rangle$ action on $\chi_i(\delta_s)$ is induced by the right $\langle c_{\lambda_i}\rangle$ action on $\delta_s$. More explicitly,  $c_{\lambda_i}$ acts on $\delta_s$  by $$(\delta_s\circ c_{\lambda_i})(j)=\left\{\begin{array}{rl}
 	\delta_s(j+1),& \text{if }\lambda_1+\cdots+\lambda_{i-1}+1\leq j \leq \lambda_1+\cdots+\lambda_i-1,\\
 	\delta_s(j-\lambda_i+1),& \text{if }j=\lambda_1+\cdots+\lambda_i,\\
 	\delta_s(j),& \text{otherwise.}
 \end{array}\right.$$

 {\noindent\bf Lemma 3.3.}
{\it Let $\lambda=\{\lambda_1\geq\lambda_2\geq\cdots\}\in\Lambda(n-*)$ and
	$$\chi_i(\delta_s)=(\underbrace{0,\cdots,0}_{k_1},a_1,\underbrace{0,\cdots,0}_{k_2},a_2,\cdots,\underbrace{0,\cdots,0}_{k_m},a_m)$$
	with $\lambda_i=2d_i$ and $a_1+\cdots+a_m=d_i$.
	$\chi_i(\delta_s)=\overline{\chi_i}(\delta_s)$ if and only if there exists a canonical
	representative of $s$ such that
	$$\delta_s(\lambda_1+\cdots+\lambda_i+1-j)=\delta_s(\lambda_1+\cdots+\lambda_{i-1}+j)+1$$ mod $2$ for every $1\leq j\leq d_i$.
}

{\noindent\bf Proof.}
By the definition of $\sigma$, we obtain $\delta_{\sigma s}(j)=\delta_s(j)+1$ mod $2$ for every $1\leq j \leq 2q$.

Assume that $$\delta_s(\lambda_1+\cdots+\lambda_i+1-j)=\delta_s(\lambda_1+\cdots+\lambda_{i-1}+j)+1$$ mod $2$ for every $1\leq j\leq d_i$, then
we have $$\delta_{\sigma s}(\lambda_1+\cdots+\lambda_i+1-j)=\delta_s((\lambda_1+\cdots+\lambda_{i-1}+j).$$ From Definition $2.2$, we consider 
$B_i(\delta_s)=\{b_1<\cdots<b_{d_i}\}$. Then we have $B_i(\delta_{\sigma s})=\{2\lambda_1+\cdots+2\lambda_{i-1}+\lambda_i+1-b_{d_i}<\cdots<2\lambda_1+\cdots+2\lambda_{i-1}+\lambda_i+1-b_1\}$. Hence we have
\begin{align*}
	\chi_i(\delta_{\sigma s})=&(b_{d_i}-b_{d_{i-1}}-1,b_{d_{i-1}}-b_{d_{i-2}}-1,\cdots,b_2-b_1-1,\lambda_i-b_{d_i}+b_1-1)\\
	=&(\lambda_i-b_{d_i}+b_1-1,b_{d_i}-b_{d_{i-1}}-1,\cdots,b_2-b_1-1)\\
	=&\overleftarrow{\chi_i(\delta_s)}.
\end{align*}
Thus we get $\chi_i(\delta_s)=\overleftarrow{\chi_i(\delta_{\sigma s})}=\sigma\circ\chi_i(\delta_s)=\overline{\chi_i}(\delta_s)$.	

Conversely, assume that $\chi_i(\delta_s)=\overline{\chi_i}(\delta_s)$. From Lemma $3.1$, as an invariant cycle,
we have
\begin{align*}\overleftarrow{\chi_i(\delta_{\sigma s})}=&(\underbrace{0,\cdots,0}_{a_m-1},,k_m+1,\cdots,\underbrace{0,\cdots,0}_{a_2-1},k_2+1,\underbrace{0,\cdots,0}_{a_1-1},k_1+1)\\
=&(\underbrace{0,\cdots,0}_{a_{m-1}-1},,k_{m-1}+1,\cdots,\underbrace{0,\cdots,0}_{a_1-1},k_1+1,\underbrace{0,\cdots,0}_{a_m-1},k_m+1)\\
=&\cdots\\
=&(\underbrace{0,\cdots,0}_{a_1-1},k_1+1,\underbrace{0,\cdots,0}_{a_m-1},k_m+1,\cdots,\underbrace{0,\cdots,0}_{a_2-1},k_2+1)\\
=&\chi_i(\delta_s)\\
=&(\underbrace{0,\cdots,0}_{k_1},a_1,\underbrace{0,\cdots,0}_{k_2},a_2,\cdots,\underbrace{0,\cdots,0}_{k_m},a_m).
\end{align*}
Hence  there exists some $1\leq i_0\leq m$, such that
$$\begin{array}{rl}
	a_{i_0-j}-1=k_{1+j},& \text{ for } 0\leq j\leq i_0-1,\\
	a_{m-j}-1=k_{i_0+j+1},& \text{ for } 0\leq j\leq m-i_0-1.
\end{array}$$
Note that 
\begin{align*}
	\chi_i(\delta_s)=&(\underbrace{0,\cdots,0}_{k_1},a_1,\underbrace{0,\cdots,0}_{k_2},a_2,\cdots,\underbrace{0,\cdots,0}_{k_m},a_m)\\
	=&(\underbrace{0,\cdots,0}_{k_{i_0+1}},a_{i_0+1},\cdots,\underbrace{0,\cdots,0}_{k_{m}},a_{m},\underbrace{0,\cdots,0}_{k_{1}},a_{1},\cdots,\underbrace{0,\cdots,0}_{k_{i_0}},a_{i_0})\\
	=&(\underbrace{0,\cdots,0}_{a_{m}-1},a_{i_0+1},\cdots,\underbrace{0,\cdots,0}_{a_{i_0+1}-1},a_{m},\underbrace{0,\cdots,0}_{a_{i_0}-1},a_{1},\cdots\underbrace{0,\cdots,0}_{a_1-1},a_{i_0}).
\end{align*}
Hence we can take a canonical representative of $s$ such that for $\lambda_1+\cdots+\lambda_{i-1}+1\leq j\leq \lambda_1+\cdots+\lambda_{i}$, $\delta_s$ satisfies
$\delta_s(j)=1$ for 
$$\begin{array}{rcl}
\lambda_1+\cdots+\lambda_{i-1}+1&\leq j \leq& \lambda_1+\cdots+\lambda_{i-1}+a_m,\\
\lambda_1+\cdots+\lambda_{i-1}+a_m+a_{i_0}+1&\leq j \leq& \lambda_1+\cdots+\lambda_{i-1}+a_m+a_{i_0}+a_{m-1},\\
&\cdots&\\
\lambda_1+\cdots+\lambda_{i}-a_1-a_{i_0}+1&\leq j \leq &\lambda_1+\cdots+\lambda_{i}-a_{i_0}
\end{array}	
$$
and $\delta_s(j)=0$ for otherwise $\lambda_1+\cdots+\lambda_{i-1}+1\leq j\leq \lambda_1+\cdots+\lambda_{i}$. Note that for such $s$, we have
$\delta_s\circ c_{\lambda_i}^{-a_1-\cdots-a_{i_0}}$ satisfies 	$$\delta_s\circ c_{\lambda_i}^{-a_1-\cdots-a_{i_0}}(\lambda_1+\cdots+\lambda_i+1-j)=\delta_s\circ c_{\lambda_i}^{-a_1-\cdots-a_{i_0}}(\lambda_1+\cdots+\lambda_{i-1}+j)+1$$ mod $2$ for every $1\leq j\leq d_i$. Hence we can take a canonical representative $s'$ of $s$ such that $\delta_{s'}=\delta_s\circ c_{\lambda_i}^{-a_1-\cdots-a_{i_0}}$.

The Lemma follows.\qb 
 
 For a fixed $\lambda=\{\lambda_,\geq\lambda_2\geq\cdots\}\in\Lambda(n-*)$, we denote
 $$\overline{\Pi(\lambda_i,\frac{\lambda_i}{2})}=\left\{\chi_i(\delta_s)\left|\begin{array}{c}(s\otimes\phi_{\lambda_i}\varepsilon,1)_{(Z_{\lambda})_s}=1,\\ \chi_i(\delta_s)=\overline{\chi_i}(\delta_s)\end{array}\right.\right\}$$
 and write simply $\overline{\Pi(2d_i,d_i)}$ for some $d_i\in \mathbb{Z}_+$.

 {\noindent\bf Lemma 3.4.} {\it Let $\lambda_i=2d_i\geq 2$ and $s\otimes\zeta_{\lambda}\in H^{\lambda_i-1}(P_n)^{\mathfrak{S}_{q} \times \mathfrak{S}_q}$ such that $$\chi(\delta_s)=\{\chi_i(\delta_s),(0),\cdots,(0),\emptyset_1,\cdots,\emptyset_1\}$$
 with $\chi_i(\delta_s)\in \overline{\Pi(2d_i,d_i)}$.	
 Then $$\sigma\circ(s\otimes\zeta_{\lambda})=-1^{\frac{(\lambda_i-1)(\lambda_i-2)}{2}} s\otimes\zeta_{\lambda}\in H^{\lambda_i-1}(P_n)^{\mathfrak{S}_{q} \times \mathfrak{S}_q}.$$}
 
 {\bf Proof:} We consider the $\sigma$ action on $s\otimes\zeta_\lambda$. We recall that from the study of Arnold \cite{A69} and by the surjection 
\[ \mathcal{A}_{\mathfrak{S}_{q} \times \mathfrak{S}_q}: H^{*}(P_n)\twoheadrightarrow
  H^{*}(P_n)^{\mathfrak{S}_{q} \times \mathfrak{S}_q},
\]
$s\otimes\zeta_\lambda$ equals the image  $$\mathcal{A}_{\mathfrak{S}_{q} \times \mathfrak{S}_q}\circ(\omega_{s(1),s(2)}\omega_{s(2),s(3)}\cdots\omega_{s(\lambda_i-1),s(\lambda_i)}).$$
Note that 
$$\sigma\circ\mathcal{A}_{\mathfrak{S}_{q} \times \mathfrak{S}_q}=\mathcal{A}_{\mathfrak{S}_{q} \times \mathfrak{S}_q}\circ\sigma.$$

 If $\lambda_i=2$, then $\chi_{i}(\delta_s)=(1)$. Hence 
 $s\otimes\zeta_\lambda$ equals the image   $\mathcal{A}_{\mathfrak{S}_{q} \times \mathfrak{S}_q}\circ\omega_{1,2q}$. We get  $$\sigma\circ(s\otimes\zeta_{\lambda})= s\otimes\zeta_{\lambda}$$
 immediately.

Now let $\lambda_i\geq3$. We write simply $\chi_i(\delta_s)\in \overline{\Pi(2d,d)}$.  From Lemma $3.3$, we take a canonical representative of $s$ such that 
$\delta_{s}(\lambda_i+1-j)=\delta_{s}(j)+1$ mod $2$ for every $1\leq j \leq d$, then
$$\omega_{s(1),s(2)}\omega_{s(2),s(3)}\cdots\omega_{s(\lambda_i-1),s(\lambda_i)}$$ is an element given by Arnold in \cite{A69} associated with $s\otimes\zeta_\lambda$.
We have
\begin{align*}
	&\sigma\circ\omega_{s(1),s(2)}\omega_{s(2),s(3)}\cdots\omega_{s(\lambda_i-1),s(\lambda_i)}\\
	=&\omega_{\sigma s(1),\sigma s(2)}
	\omega_{\sigma s(2),\sigma s(3)}
	\cdots
	\omega_{\sigma s(\lambda_i-1),\sigma s(\lambda_i)}\\
	=&-1^{\frac{(\lambda_i-1)(\lambda_i-2)}{2}}
	\omega_{\sigma s(\lambda_i),\sigma s(\lambda_i-1)}
	\cdots
	\omega_{\sigma s(3),\sigma s(2)}
	\omega_{\sigma s(2),\sigma s(1)}.
\end{align*}

Note that $\delta_{\sigma s}(\lambda_i+1-j)=\delta_{s}(j)$ for every $1\leq j \leq d$,	
  Hence 
   there exists a $\tau\in\mathfrak{S_q}\times \mathfrak{S_q}$, such that 
\begin{align*}
 &\omega_{\sigma s(\lambda_i),\sigma s(\lambda_i-1)}\cdots\omega_{\sigma s(3),\sigma s(2)}\omega_{\sigma s(2),\sigma s(1)}\\
 =&\tau\circ\omega_{ s(1), s(2)}\omega_{s(2),s(3)}\cdots\omega_{ s(\lambda_i-1), s(\lambda_i)}.
\end{align*}
    Note that
\begin{align*}    
&\mathcal{A}_{\mathfrak{S}_{q} \times \mathfrak{S}_q}\circ(\omega_{ s(1), s(2)}\omega_{s(2), s(3)}\cdots\omega_{ s(\lambda_i-1), s(\lambda_i)})\\
=&\mathcal{A}_{\mathfrak{S}_{q} \times \mathfrak{S}_q}\circ(\tau\circ\omega_{ s(1), s(2)}\omega_{ s(2), s(3)}\cdots\omega_{ s(\lambda_i-1), s(\lambda_i)})\\
=&\mathcal{A}_{\mathfrak{S}_{q} \times \mathfrak{S}_q}\circ(\omega_{\sigma s(\lambda_i),\sigma s(\lambda_i-1)}\cdots\omega_{\sigma s(3),\sigma s(2)}\omega_{\sigma s(2),\sigma s(1)}).
\end{align*}
 Hence we get
 $$\sigma\circ(s\otimes\zeta_{\lambda})=-1^{\frac{(\lambda_i-1)(\lambda_i-2)}{2}} s\otimes\zeta_{\lambda}.$$ 
 
 These prove the Lemma.\qb
 
 We regard $s\otimes\zeta_{\lambda}$ as an element in $\bigotimes_{i=1}^{a}H^{\lambda_i-1}(P_n)^{\mathfrak{S}_{q} \times \mathfrak{S}_q}$,  then from Lemma $3.1$ and Definition $3.2$, we can get $$\overline{\chi}(\delta_s)=\{\overline{\chi_1}(\delta_s),\overline{\chi_2}(\delta_s),\cdots,\overline{\chi_a}(\delta_s),(0),\cdots,(0),\emptyset_1,\cdots,\emptyset_1\}$$
  by rearranging the sequence 
$$(\overline{\chi_1}(\delta_s),\overline{\chi_2}(\delta_s),\cdots,\overline{\chi_a}(\delta_s))$$ by the lexicographic order we defined in Definition $2.2$. We denote by $\overline{s}$ the representative associated with $\overline{\chi}(\delta_s)$. Hence by applying Lemma $3.3$, we have 
\begin{align*}
\sigma\circ(s\otimes\zeta_{\lambda})
=\varepsilon(\lambda,s)(\overline{s}\otimes\zeta_{\lambda}),
\end{align*}
where $\varepsilon(\lambda,s)=\pm 1$ is depending on the sign of the rearrangement of the sequence $$(\overline{\chi_1}(\delta_s),\overline{\chi_2}(\delta_s),\cdots,\overline{\chi_a}(\delta_s))$$ and the sign of the $\sigma$ action on each element in $H^{\lambda_i-1}(P_n)^{\mathfrak{S}_{q} \times \mathfrak{S}_q}$.

We take $$\bigcup_{\lambda}\{s\otimes\zeta_{\lambda}|(s\otimes\zeta_{\lambda},1)_{(Z_{\lambda})_s}=1\}$$ as a basis of $H^{*}(P_n)^{\mathfrak{S}_{q} \times \mathfrak{S}_q}$.
Recall that $$\phi(s\otimes\zeta_{\lambda})=s\otimes\zeta_{\lambda}+\varepsilon(\lambda,s)(\overline{s}\otimes\zeta_{\lambda}).$$
 For a fixed  $s\otimes\zeta_{\lambda}$ of $H^{*}(P_n)^{\mathfrak{S}_{q} \times \mathfrak{S}_q}$ with a fixed partition $\lambda$,  if $\overline{s}\neq s$ or $\overline{s}=s$ with $\varepsilon(\lambda,s)=1$, then $\phi(s\otimes\zeta_{\lambda})\neq0$. If $\overline{s}\neq s$, then $\phi(\overline{s}\otimes\zeta_{\lambda})=\pm\phi(s\otimes\zeta_{\lambda})$.
 
  We note that $\overline{s}$ is uniquely determined by $s$. For a fixed $\lambda\in \Lambda(n-*)$,  we take canonical representatives of $s$ such that $\chi(\delta_s)>\overline{\chi}(\delta_s)$ for those $s\neq\overline{s} $, then the condition $\chi(\delta_s)>\overline{\chi}(\delta_s)$ selects exactly one representative from each  $\sigma$-orbit of the $\sigma$ action on $s\otimes\zeta_{\lambda}$ such that $s\neq\overline{s} $. Hence  the union of
 $$\mathfrak{B}_1=\bigcup_{\lambda\in\Lambda(n-*)}\{\phi(s\otimes\zeta_{\lambda})|(s\otimes\zeta_{\lambda},1)_{(Z_{\lambda})_s}=1, \chi(\delta_s)>\overline{\chi}(\delta_s)\}$$
 and
 $$\mathfrak{B}_2=\bigcup_{\lambda\in\Lambda(n-*)}\{\phi(s\otimes\zeta_{\lambda})|(s\otimes\zeta_{\lambda},1)_{(Z_{\lambda})_s}=1, \chi(\delta_s)=\overline{\chi}(\delta_s), \varepsilon(\lambda,s)=1\}$$
  forms a basis of $H^{*}(P_n)^{\mathfrak{G}}$. The kernel of $\phi$  is spanned by
 $$\bigcup_{\lambda}\{s\otimes\zeta_{\lambda}-\varepsilon(\lambda,s)(\overline{s}\otimes\zeta_{\lambda})|(s\otimes\zeta_{\lambda},1)_{(Z_{\lambda})_s}=1\}.$$

Now we consider those $s\otimes \zeta_{\lambda}=\overline{s}\otimes \zeta_{\lambda}$.
For a fixed partition
$$\lambda=(\lambda_1=\lambda_2\cdots=\lambda_{|(\lambda_1)|}>\lambda_{|(\lambda_1)|+1}=\lambda_{|(\lambda_1)|+2}=\cdots>1=\cdots=1)\in \Lambda(n-*),$$ we write simply
 $$\lambda=(\underbrace{\lambda_1,\cdots,\lambda_1}_{|(\lambda_1)|},\underbrace{\lambda_2,\cdots,\lambda_2}_{|(\lambda_2)|},\cdots,\underbrace{\lambda_a,\cdots,\lambda_a}_{|(\lambda_a)|},1,\cdots,1)\in \Lambda(n-*),$$
with
$$\lambda_i\in \{\lambda_1>\lambda_2>\cdots>\lambda_a\geq2\}$$
and
$$|(\lambda_i)|=\sharp\{1\leq j\leq |\lambda_1|+\cdots+|\lambda_a||\lambda_j=\lambda_{i}\}.$$

 In this section, we write simply $$l_i=\sum_{j=1}^{i-1}|(\lambda_j)|.$$

{\noindent\bf Definition 3.5.} {\it Let $$\lambda=(\underbrace{\lambda_1,\cdots,\lambda_1}_{|(\lambda_1)|},\underbrace{\lambda_2,\cdots,\lambda_2}_{|(\lambda_2)|},\cdots,\underbrace{\lambda_a,\cdots,\lambda_a}_{|(\lambda_a)|},1,\cdots,1)\in \Lambda(n-*).$$
For a $$\lambda_i\in \{\lambda_1,\lambda_2,\cdots,\lambda_a\}$$ and a $ 0\leq k_i \leq\lfloor\frac{|(\lambda_i)|}{2}\rfloor$,
if there exists some 
\begin{align*}
&((d_{l_i+1},\cdots,d_{l_i+k_i}) ,(d_{l_i+k_i+1}, \cdots, d_{l_{i+1}-k_i}),(d_{l_{i+1}-k_i+1},\cdots,d_{l_{i+1}}))\\
\in&(\mathbb{Z}_{\geq0})^{k_i}\times(\mathbb{Z}_{\geq0})^{|(\lambda_i)|-2k_i}\times(\mathbb{Z}_{\geq0})^{k_i}\cong(\mathbb{Z}_{\geq0})^{|(\lambda_i)|}
\end{align*}
such that 
$$d_{l_i+j}=\lambda_i-d_{l_{i+1}-j+1},$$ 
$$d_{l_i+j}\geq \frac{\lambda_i}{2},$$ for $1\leq j \leq k_i$ and
$$d_{l_i+j}=\frac{\lambda_i}{2}$$
for $k_i+1\leq j\leq |(\lambda_{i})|-k_i$,
then we define an element $$d_{\lambda_i}(k_i)=(k_i,(d_{l_i+1},\cdots,d_{l_i+k_i}) ,(d_{l_i+k_i+1}, \cdots, d_{l_{i+1}-k_i}),(d_{l_{i+1}-k_i+1},\cdots,d_{l_{i+1}}))$$ in $\{k_i\}\times(\mathbb{Z}_{\geq0})^{k_i}\times(\mathbb{Z}_{\geq0})^{|(\lambda_i)|-2k_i}\times(\mathbb{Z}_{\geq0})^{k_i}$ and  denote
 $D_{\lambda_i}(k_i)$ the set of all $d_{\lambda_i}(k_i)$.}

{\noindent\bf Definition 3.6.} {\it Let $$\lambda=(\underbrace{\lambda_1,\cdots,\lambda_1}_{|(\lambda_1)|},\underbrace{\lambda_2,\cdots,\lambda_2}_{|(\lambda_2)|},\cdots,\underbrace{\lambda_a,\cdots,\lambda_a}_{|(\lambda_a)|},1,\cdots,1)\in \Lambda(n-*).$$ If there exists some $ 0\leq k_i \leq\lfloor\frac{|(\lambda_i)|}{2}\rfloor$ such that
	$D_{\lambda_i}(k_i)\neq\emptyset$ for a $$\lambda_i\in \{\lambda_1,\lambda_2,\cdots,\lambda_a\},$$
then we define a multiset of length $|(\lambda_i)|$, denoted by $(\lambda_i,d_{\lambda_i}(k_i))$, by 
$$(\lambda_i,d_{\lambda_i}(k_i))=\{(\lambda_i,k_i,d_{l_i+1}),\cdots, (\lambda_i,k_i,d_{l_{i+1}})\}$$  
and define a set of multisets, denoted by $E_{\lambda_i}$, by
$$E_{\lambda_i}=\bigsqcup_{k_i=0}^{  \lfloor\frac{|(\lambda_i)|}{2}\rfloor}\{(\lambda_i,d_{\lambda_i}(k_i))|d_{\lambda_i}(k_i)\in D_{\lambda_i}(k_i)\neq\emptyset\}.$$
If 
$E_{\lambda_i}\neq\emptyset$ for every $$\lambda_i\in \{\lambda_1,\lambda_2,\cdots,\lambda_a\},$$
then we also define a set of multisets of length $|(\lambda_1)|+\cdots+|(\lambda_a)|$, denoted by $E_\lambda$,
by $$E_\lambda=\left\{\bigcup_{\lambda_i\in\{\lambda_1,\cdots,\lambda_a\}}(\lambda_i,d_{\lambda_i}(k_i))\left|
\begin{array}{c}
\lambda\in\Lambda(n-*),d_1+\cdots+d_{|(\lambda_1)|+\cdots+|(\lambda_a)|}\leq q,\\
(\lambda_i,d_{\lambda_i}(k_i))\in E_{\lambda_i}	
\end{array}\right.
\right\}.$$
We write simply $$(\lambda,d,k)= \bigcup_{\lambda_i\in\{\lambda_1,\cdots,\lambda_a\}}(\lambda_i,d_{\lambda_i}(k_i))\in E_\lambda$$ for $k=(k_1,\cdots,k_a)$.
}

Recall that $$(\lambda,d)=\{(\lambda_1,d_1),\cdots,(\lambda_{|(\lambda_1)|+\cdots+|(\lambda_a)|},d_{|(\lambda_1)|+\cdots+|(\lambda_a)|})\}$$ is a multiset of length $|(\lambda_1)|+\cdots+|(\lambda_a)|$. If $E_\lambda\neq\emptyset$,  then there exists some  $(\lambda,d,k)\in E_\lambda$ such that  $(\lambda,d,k)$ is associated with $(\lambda,d)$.
 
{\noindent\bf Definition 3.7.} {\it Let $$\lambda=(\underbrace{\lambda_1,\cdots,\lambda_1}_{|(\lambda_1)|},\underbrace{\lambda_2,\cdots,\lambda_2}_{|(\lambda_2)|},\cdots,\underbrace{\lambda_a,\cdots,\lambda_a}_{|(\lambda_a)|},1,\cdots,1)\in\Lambda(n-*)$$ and $E_{\lambda}\neq\emptyset$.
	 We define a subset of the basis of $H^{*}(P_n)^{\mathfrak{S}_q\times\mathfrak{S}_q}$,
	denoted by $E_{\lambda}^P$,
	 by  \begin{align*}&E_{\lambda}^P
	=&\bigcup_{(\lambda,d,k)\in E_\lambda}\left\{s\otimes\zeta_{\lambda}\left|\begin{array}{c}
		s\otimes\zeta_{\lambda} \text{ corresponds to the fixed } (\lambda,d,k),\\
		(s\otimes\zeta_{\lambda},1)_{(Z_{\lambda})_s}=1,\\
		\exists \tau \in \mathfrak{S}_{|(\lambda_1)|}\times\cdots\times\mathfrak{S}_{|(\lambda_a)|}\subset \mathfrak{S}_{n}, \text{ s.t.}\\
		\overline{\chi_{\tau(l_i+j)}}(\delta_s)=\chi_{\tau(l_{i+1}-j+1)}(\delta_s),\\
		\chi_{\tau(l_i+j)}(\delta_s)<\overline{\chi_{\tau(l_i+j)}}(\delta_s),\\
		1\leq i\leq a, 1\leq j\leq k_i,\\
		\chi_{\tau(l_i+j)}(\delta_s)\leq\chi_{\tau(l_i+j+1)}(\delta_s),\\
		1\leq i\leq a, 1\leq  j\leq k_i-1,\\
		\chi_{\tau(l_i+j)}(\delta_s)=\overline{\chi_{\tau(l_i+j)}}(\delta_s),\\
			1\leq i\leq a,  k_i+1\leq j \leq |(\lambda_{i})|-k_i,\\
		\chi_{\tau(l_i+j)}(\delta_s)<\chi_{\tau(l_i+j+1)}(\delta_s),\\
		1\leq i\leq a,  k_i+1\leq j \leq |(\lambda_{i})|-k_i-1
	\end{array}\right.\right\},\end{align*}
where $\chi(\delta_s)$ follows the lexicographic order defined in Definition $2.2$ and $\tau$ is a relabeling of $\chi(\delta_s)$ fixed by $(\lambda,d,k)\in E_\lambda$ such that $\lambda_{\tau(j)}=\lambda_{j}$ and $d_{\tau(j)}=d_j$, $1\leq j \leq l_{a+1}$.}

{ \bf Remark:} From Lemma $3.1$, for a fixed $(\lambda,d,k)\in E_\lambda$, $s\otimes \zeta_{\lambda}\in E_\lambda^P$ is uniquely determined by $$\chi_{\tau(l_i+j)}(\delta_s),$$ $1\leq \tau(j) \leq |(\lambda_{i})|-k_i$ for each $\lambda_i$, $1\leq i\leq a$. 

From now we relabel $\chi(\delta_s)$ for every $s\otimes\zeta_{\lambda}\in E_{\lambda}^P$ by $\chi_{j}(\delta_s)=\chi_{\tau(j)}(\delta_s)$, $1\leq j \leq l_{a+1}$, i.e., we write
\begin{align*}
&\overline{\chi_{l_i+j}}(\delta_s)=\chi_{l_{i+1}-j+1}(\delta_s),
\chi_{l_i+j}(\delta_s)<\overline{\chi_{l_i+j}}(\delta_s),
1\leq i\leq a, 1\leq j\leq k_i,\\
&\chi_{l_i+j}(\delta_s)\leq\chi_{l_i+j+1}(\delta_s),
1\leq i\leq a, 1\leq j\leq k_i-1,\\
&\chi_{l_i+j}(\delta_s)=\overline{\chi_{l_i+j}}(\delta_s),
1\leq i\leq a,  k_i+1\leq j \leq |(\lambda_{i})|-k_i,\\
&\chi_{l_i+j}(\delta_s)<\chi_{l_i+j+1}(\delta_s),
1\leq i\leq a,  k_i+1\leq j \leq |(\lambda_{i})|-k_i-1
\end{align*}
if $s\otimes\zeta_{\lambda}\in E_{\lambda}^P$. 
The condition  may not follow the lexicographic order  defined in Definition $2.2$, but the relabeling does not change the underlying full invariant cycle set. Since all factors in a fixed $\lambda_i$ have varvathe same degree $\lambda_i-1$, the relabeling conjugates the permutation induced by $\sigma$. Hence  it does not change the parity of that 
permutation or the resulting sign $\varepsilon(\lambda,s)$. We will use $E_\lambda^P$ with this condition to make the computations of 
 $\overline{\chi}(\delta_s)=\chi(\delta_s)$ and $\varepsilon(\lambda,s)=-1$ easier.

{\noindent \bf Lemma 3.8.} {\it Let $s\otimes \zeta_{\lambda}$ be a basis element of $H^*(P_n)^{\mathfrak{S}_q\times\mathfrak{S}_q}$. $s\otimes\zeta_{\lambda}=\overline{s}\otimes\zeta_{\lambda}$ if and only if $s\otimes\zeta_{\lambda}\in E_\lambda^P$.}

{\noindent\bf Proof.} Assume that $s\otimes\zeta_{\lambda}\in E_\lambda^P$, we consider the $\sigma$ action on $\chi(\delta_s)$ associated with $s\otimes\zeta_{\lambda}$, i.e., we consider  $\overline{\chi}(\delta_s)$.
We have
$$
 \overline{\overline{\chi_{l_i+j}}}(\delta_{s})=\chi_{l_i+j}(\delta_{s})
$$
for $\lambda_i\in\{\lambda_1,\cdots,\lambda_a\}$, $1\leq j\leq |(\lambda_i)|$.
  By  definition $3.7$, we have
$$\begin{array}{rl}
	&\overline{\overline{\chi_{l_{i+1}-j+1}}}(\delta_{ s})=\chi_{l_{i+1}-j+1}(\delta_s)=\overline{\chi_{l_i+j}}(\delta_s)\\
	>&\chi_{l_i+j}(\delta_s)=\overline{\overline{\chi_{l_{i}+j}}}(\delta_{ s})=\overline{\chi_{l_{i+1}-j+1}}(\delta_{ s})
\end{array}$$
for $1\leq i\leq a$, $1\leq j\leq k_i$.
$$\begin{array}{rl}
	&\overline{\chi_{l_{i+1}-j+1}}(\delta_{s})=\overline{\overline{\chi_{l_i+j}}}(\delta_s)\\
	\leq&\overline{\overline{\chi_{l_i+j+1}}}(\delta_s)=\overline{\chi_{l_{i+1}-j}}(\delta_{s})
\end{array}$$ for $1\leq i\leq a$, $1\leq j\leq k_i-1$.
$$\begin{array}{rl}
	&\overline{\overline{\chi_{l_i+j}}}(\delta_s)=\overline{\chi_{l_i+j}}(\delta_{ s})\\
	<&\overline{\overline{\chi_{l_i+j+1}}}(\delta_s)=\overline{\chi_{l_i+j+1}}(\delta_{s})
\end{array}$$
for $1\leq i\leq a$,  $k_i+1\leq j \leq  |(\lambda_{i})|-k_i-1$.
Note that $\overline{\chi}(\delta_s)$ and $\chi(\delta_s)$ are multisets, hence we have  $\overline{\chi}(\delta_s)=\chi(\delta_s)$. 
 Thus as an element of $H^*(P_n)^{\mathfrak{S}_q\times\mathfrak{S}_q}$, we have $s\otimes\zeta_{\lambda}=\overline{s}\otimes\zeta_{\lambda}$.

Conversely, if $s\otimes\zeta_{\lambda}=\overline{s}\otimes\zeta_{\lambda}$, then the full invariant cycle set associated with $\sigma \circ (s\otimes \zeta_{\lambda})=\pm\overline{s}\otimes\zeta_{\lambda}$ satisfies $\overline{\chi}(\delta_s)=\chi(\delta_s)$. Hence there exists some $ 0\leq k_i \leq\lfloor\frac{|(\lambda_i)|}{2}\rfloor$, $1\leq i \leq a$, such that there exists $k_i$ invariant cycles $\chi_{r_l}(\delta_s)$ satisfy $\overline{\chi_{r_t}}(\delta_s)=\chi_{r_l}(\delta_s)$ and there exists $|(\lambda_i)|-2k_i$  invariant cycles $\chi_{r_l}(\delta_s)$ satisfy $\overline{\chi_{r_l}}(\delta_s)=\chi_{r_l}(\delta_s)$, where $r_l\neq r_t$ is fixed by $\lambda_i$. By relabeling such $\chi(\delta_s)$, we have $s\otimes\zeta_{\lambda}\in E_\lambda^P$.
\qb 

From Lemma $3.8$,  $E_\lambda^P$ constitutes a subset of the  basis of $H^*(P_n)^{\mathfrak{S}_q\times\mathfrak{S}_q}$ for a fixed $\lambda$.
Our subsequent computations of $H^*(P_n)^{\mathfrak{G}}$ will rely on $E_\lambda^P$.

{\noindent \bf Lemma 3.9.} {\it Let  $\phi:H^{*}(P_n)^{\mathfrak{S}_{q} \times \mathfrak{S}_q}\twoheadrightarrow	H^{*}(P_n)^{\mathfrak{G}}$ be the surjection such that $$\phi(s\otimes\zeta_{\lambda})=s\otimes\zeta_{\lambda}+\varepsilon(\lambda,s)(\overline{s}\otimes\zeta_{\lambda}).$$ If $s\otimes\zeta_{\lambda}\in E_\lambda^P$, then $$\varepsilon(\lambda,s)=\prod_{\lambda_i\in \{\lambda_1,\lambda_2,\cdots,\lambda_a\}}(-1)^{\frac{(|(\lambda_i)|-2k_i)(\lambda_i-1)(\lambda_i-2)}{2}+(\lambda_i-1)k_i(2k_i-1)}.$$}

{\noindent\bf Proof.}
 For a fixed $$\lambda=(\underbrace{\lambda_1,\cdots,\lambda_1}_{|(\lambda_1)|},\underbrace{\lambda_2,\cdots,\lambda_2}_{|(\lambda_2)|},\cdots,\underbrace{\lambda_a,\cdots,\lambda_a}_{|(\lambda_a)|},1,\cdots,1)\in\Lambda(n-*),$$ we assume that $s\otimes\zeta_{\lambda}\in E_\lambda^P$. 
 
 Consider a $$\lambda_i\in \{\lambda_1,\lambda_2,\cdots,\lambda_a\},$$ from Lemma $3.4$, we have the sign of the $\sigma$ action on the element associated with $$\chi_{l_i+j}(\delta_s)=\overline{\chi_{l_i+j}}(\delta_s),$$ where $j$ is fixed by $E_\lambda^P$, is $$-1^{\frac{(\lambda_i-1)(\lambda_i-2)}{2}}.$$ By  Definition $3.7$, totally we have $|(\lambda_i)|-2k_i$ such elements, hence we get
 $$-1^{\frac{(|(\lambda_i)|-2k_i)(\lambda_i-1)(\lambda_i-2)}{2}}.$$
 
It is notable that we don't need to consider the sign of the $\sigma$ action on those elements in $H^{\lambda_i-1}(P_n)^{\mathfrak{S}_q\times\mathfrak{S}_q}$ such that $\chi_{l_i+j}(\delta_s)\neq\overline{\chi_{l_i+j}}(\delta_s)$.
 
Now we consider the sign of the rearrangement of the sequence
$$(\overline{\chi_{l_i+1}}(\delta_s), \cdots,\overline{\chi_{l_{i+1}}}(\delta_s)).$$ Recall that we have 
$$\begin{array}{rl}
	&\overline{\overline{\chi_{l_{i+1}-j+1}}}(\delta_{ s})=\chi_{l_{i+1}-j+1}(\delta_s)=\overline{\chi_{l_i+j}}(\delta_s)\\
	>&\chi_{l_i+j}(\delta_s)=\overline{\overline{\chi_{l_{i}+j}}}(\delta_{ s})=\overline{\chi_{l_{i+1}-j+1}}(\delta_{ s})
\end{array}$$
for $1\leq i\leq a$, $1\leq j\leq k_i$.
$$\begin{array}{rl}
	&\overline{\chi_{l_{i+1}-j+1}}(\delta_{s})=\overline{\overline{\chi_{l_i+j}}}(\delta_s)\\
	\leq&\overline{\overline{\chi_{l_i+j+1}}}(\delta_s)=\overline{\chi_{l_{i+1}-j}}(\delta_{s})
\end{array}$$for $1\leq i\leq a$, $1\leq j\leq k_i-1$.
Then we can get the sign by rearranging $$(\overline{\chi_{l_i+1}}(\delta_s), \cdots,\overline{\chi_{l_{i+1}}}(\delta_s))$$ to 
$$(\chi_{l_i+1}(\delta_s),\cdots,\chi_{l_{i+1}}(\delta_s)).$$
We note that
$$\begin{array}{rl}
(\chi_{l_i+1}(\delta_s),\cdots,\chi_{l_{i}+k_i}(\delta_s))
&=(\overline{\chi_{l_{i+1}}}(\delta_s),\cdots,\overline{\chi_{l_{i+1}-k_i+1}}(\delta_s)),\\

(\chi_{l_{i}+k_i+1} \cdots, \chi_{l_{i+1}-k_i}(\delta_s))
&=(\overline{\chi_{l_i+k_i+1}}(\delta_s),\cdots,\overline{\chi_{l_{i+1}-k_i}}(\delta_s)),\\

(\chi_{l_{i+1}-k_i+1} \cdots, \chi_{l_i+1}(\delta_s))
&=(\overline{\chi_{l_i+k_i}}(\delta_s),\cdots,\overline{\chi_{l_{i}+1}}(\delta_s)).
\end{array}$$
Hence the sign of this rearrangement is
$$-1^{(\lambda_i-1)(\frac{2k_i(2k_i-1)}{2}+2(|(\lambda_i)|-2k_i)k_i)}=-1^{(\lambda_i-1)\frac{2k_i(2k_i-1)}{2}}=-1^{(\lambda_i-1)k_i(2k_i-1)}.$$

By running over all $$\lambda_i\in\{\lambda_1,\cdots,\lambda_a\},$$ we get our Lemma. \qb

Now consider the kernel of the surjection $$\phi: H^{*}(P_n)^{\mathfrak{S}_{q} \times \mathfrak{S}_q}\rightarrow H^{*}(P_n)^{\mathfrak{G}}.$$ Recall that The kernel of $\phi$  is spanned by
$$\bigcup_{\lambda}\{s\otimes\zeta_{\lambda}-\varepsilon(\lambda,s)(\overline{s}\otimes\zeta_{\lambda})|(s\otimes\zeta_{\lambda},1)_{(Z_{\lambda})_s}=1\}.$$
For a fixed $\lambda\in\Lambda(n-*)$, if 
$$s\otimes\zeta_{\lambda}-\varepsilon(\lambda,s)(\overline{s}\otimes\zeta_{\lambda})=2s\otimes\zeta_{\lambda},$$
then from Lemma $3.8$, we have $s\otimes\zeta_{\lambda}\in E_\lambda^P$. We write 
$$E^P=\bigcup_{\lambda\in\Lambda(n-*)}E_\lambda^P.$$
We only consider the kernel of $\phi$ restricted by $E^P$.

{\noindent\bf Theorem 3.10.} {\it Let	
	  $$K^P=\bigcup_{\lambda\in\Lambda(n-*)}\bigcup_{(\lambda,d,k)\in K_\lambda}\{s\otimes\zeta_{\lambda}\in E_\lambda^P| s\otimes\zeta_{\lambda} \text{ correspond to } (\lambda,d,k)\in K_\lambda\}\subset E^P,$$
 where $K_\lambda\subset E_\lambda$ is the set of all $(\lambda,d,k)$ such that :  there exists odd number of $1\leq i \leq a$, such that $|(\lambda_i)|$ and $k_i$ satisfy one of: 
\begin{enumerate}
	\item  If $\lambda_i\equiv0$ mod $4$, then $|(\lambda_i)|\equiv 0$ mod $2$ with $k_i\equiv1$ mod $2$ or $|(\lambda_i)|\equiv 1$ mod $2$ with $k_i\equiv0$ mod $2$. 
	\item If $\lambda_i\equiv2$ mod $4$, then $k_i\equiv1$ mod $2$.
\end{enumerate} 
Then $K^P$ is a  basis of the kernel of $\phi$ restricted by the $\mathbb{Q}$-vector space spanned by $E^P$.}

{\noindent\bf Proof.} Assume that $s\otimes\zeta_{\lambda}\in E^P$. From Lemma $3.9$, we obtain  $$\varepsilon(\lambda,s)=\prod_{\lambda_i\in \{\lambda_1,\lambda_2,\cdots,\lambda_a\}}(-1)^{\frac{(|(\lambda_i)|-2k_i)(\lambda_i-1)(\lambda_i-2)}{2}+(\lambda_i-1)k_i(2k_i-1)}.$$
We establish $\varepsilon(\lambda,s)=-1$ through the following computations.

For a fixed $$\lambda=(\underbrace{\lambda_1,\cdots,\lambda_1}_{|(\lambda_1)|},\underbrace{\lambda_2,\cdots,\lambda_2}_{|(\lambda_2)|},\cdots,\underbrace{\lambda_a,\cdots,\lambda_a}_{|(\lambda_a)|},1,\cdots,1)\in\Lambda(n-*),$$
the condition $E_\lambda^P\neq\emptyset$ implies that if $\lambda_i\equiv1$ mod $2$, then $|(\lambda_i)|\equiv0$ mod $2$ and $k_i=\frac{|(\lambda_i)|}{2}$. 
Hence for this case, $$(-1)^{\frac{(|(\lambda_i)|-2k_i)(\lambda_i-1)(\lambda_i-2)}{2}+(\lambda_i-1)k_i(2k_i-1)}=1.$$

If 
$\lambda_i\equiv0$ mod $2$, we have $$(-1)^{\frac{(|(\lambda_i)|-2k_i)(\lambda_i-1)(\lambda_i-2)}{2}+(\lambda_i-1)k_i(2k_i-1)}=-1$$ if and only if $|(\lambda_i)|$ and $k_i$ satisfy one of
\begin{enumerate}
\item  If $\lambda_i\equiv0$ mod $4$, then $|(\lambda_i)|\equiv 0$ mod $2$ with $k_i\equiv1$ mod $2$ or $|(\lambda_i)|\equiv 1$ mod $2$ with $k_i\equiv0$ mod $2$. 
\item If $\lambda_i\equiv2$ mod $4$, then $k_i\equiv1$ mod $2$.
\end{enumerate}

These prove the theorem. \qb

\section{the dimension of $H^{*}(P_n)^\mathfrak{G}$}

Recall that $\mathfrak{G}=\mathfrak{S}_q\overleftrightarrow{\times}\mathfrak{S}_q$ and $n=2q$.
In this section, we will give a formula to compute the dimension of $H^{*}(P_n)^\mathfrak{G}$.
  
{\noindent\bf Definition 4.1.}
{\it Let $$\lambda=(\underbrace{\lambda_1,\cdots,\lambda_1}_{|(\lambda_1)|},\underbrace{\lambda_2,\cdots,\lambda_2}_{|(\lambda_2)|},\cdots,\underbrace{\lambda_a,\cdots,\lambda_a}_{|(\lambda_a)|},1,\cdots,1)\in\Lambda(n-*),$$ 
$$E_\lambda=\left\{\bigcup_{\lambda_i\in\{\lambda_1,\cdots,\lambda_a\}}(\lambda_i,d_{\lambda_i}(k_i))\left|
\begin{array}{c}
	\lambda\in\Lambda(n-*),d_1+\cdots+d_{|(\lambda_1)|+\cdots+|(\lambda_a)|}\leq q,\\
	(\lambda_i,d_{\lambda_i}(k_i))\in E_{\lambda_i}	
\end{array}\right.
\right\}$$
and $$l_i=\sum_{j=1}^{i-1}|(\lambda_j)|.$$
For a fixed 
$$(\lambda,d,k)=\bigcup_{\lambda_i\in\{\lambda_1,\cdots,\lambda_a\}}(\lambda_i,d_{\lambda_i}(k_i))\in E_\lambda,$$
we define
\begin{align*}
	E_1(\lambda,d)&=\{(\lambda_j,d_j)|1\leq j \leq l_{a+1},  \lambda_j\equiv 1 \text{ mod } 2, d_j>\frac{\lambda_i}{2}, i \text{ is fixed by }j\},\\
	E_{21}(\lambda,d,k)&=\{(\lambda_j,d_j)|l_i+1\leq j \leq l_i+k_i,\lambda_j \equiv 0 \text{ mod } 2, d_j > \frac{\lambda_i}{2}, i \text{ is fixed by }j\},\\
	E_{22}(\lambda,d,k)&=\{(\lambda_j,d_j)|l_i+1\leq j \leq l_i+k_i,\lambda_j \equiv 0 \text{ mod } 2, d_j = \frac{\lambda_i}{2}, i \text{ is fixed by }j\},\\
	E_{3}(\lambda,d,k)&=\left\{(\lambda_j,d_j)\left|\begin{array}{c}l_i+k_i+1\leq j \leq l_{i+1}-k_i,\lambda_i \equiv 0 \text{ mod } 2, \\
		d_j= \frac{\lambda_i}{2},i \text{ is fixed by }j\end{array}\right.\right\}
\end{align*}
and for a fixed $(\lambda_i,d_i)$ in one of $E_1(\lambda,d)$, $E_{21}(\lambda,d,k)$, $E_{22}(\lambda,d,k)$ or $E_{3}(\lambda,d,k)$,  we define
$$|(\lambda,d)_i|=\sharp\{j|(\lambda_j,d_j)=(\lambda_{i},d_i)\},$$
where $(\lambda_j,d_j)$  are fixed by one of $E_1(\lambda,d)$, $E_{21}(\lambda,d,k)$, $E_{22}(\lambda,d,k)$ or  $E_{3}(\lambda,d,k)$.
For $1\leq b \leq |(\lambda,d)_i|$, we also denote simply
$$
M(b)=\{(m_1,m_2,\cdots,m_b)\in (\mathbb{Z}_{+})^b| m_1+m_2+\cdots+m_b=	|(\lambda,d)_i|\}
$$
for $M(\lambda,d,,i,b)$.
}  
 
We recall that
$$\overline{\Pi(\lambda_i,\frac{\lambda_i}{2})}=\left\{\chi_i(\delta_s)\left|\begin{array}{c}(s\otimes\phi_{\lambda_i}\varepsilon,1)_{(Z_{\lambda})_s}=1,\\ \chi_i(\delta_s)=\overline{\chi_i}(\delta_s)\end{array}\right.\right\}$$
and write simply $\overline{\Pi(2d,d)}$ for some $d\in \mathbb{Z}_+$.

{\noindent\bf Theorem 4.2.} {\it  Let $|E^P|$ and $|K^P|$ denote the cardinality of $E^P$ and $K^P$ respectively. Then $$\text{dim}(H^{*}(P_n)^{\mathfrak{G}})=\frac{\text{dim}(H^{*}(P_n)^{\mathfrak{S}_{q} \times \mathfrak{S}_q})+|E^P|}{2}-|K^P|.$$
}

{\noindent\bf Proof.} Recall that the union of
$$\mathfrak{B}_1=\bigcup_{\lambda\in\Lambda(n-*)}\{\phi(s\otimes\zeta_{\lambda})|(s\otimes\zeta_{\lambda},1)_{(Z_{\lambda})_s}=1, \chi(\delta_s)>\overline{\chi}(\delta_s)\}$$
and
$$\mathfrak{B}_2=\bigcup_{\lambda\in\Lambda(n-*)}\{\phi(s\otimes\zeta_{\lambda})|(s\otimes\zeta_{\lambda},1)_{(Z_{\lambda})_s}=1, \chi(\delta_s)=\overline{\chi}(\delta_s), \varepsilon(\lambda,s)=1\}$$
forms a basis of $H^{*}(P_n)^{\mathfrak{G}}$. $\mathfrak{B}_1\cap\mathfrak{B}_2=\emptyset$.

Since $\phi(s\otimes\zeta_{\lambda})\in\mathfrak{B}_1$ is the image of  $s\otimes\zeta_{\lambda}\notin E^P$ and $\phi(s\otimes\zeta_{\lambda})=\pm\phi(\overline{s}\otimes\zeta_{\lambda})$, the cardinality of $\mathfrak{B}_1$ is
$$\frac{\text{dim}(H^{*}(P_n)^{\mathfrak{S}_{q} \times \mathfrak{S}_q})-|E^P|}{2}.$$

On the other hand, Theorem $3.10$ implies that  $$E^P\setminus K^P=\bigcup_{\lambda\in\Lambda(n-*)}\{s\otimes\zeta_{\lambda}|(s\otimes\zeta_{\lambda},1)_{(Z_{\lambda})_s}=1, \chi(\delta_s)=\overline{\chi}(\delta_s), \varepsilon(\lambda,s)=1\}.$$ Hence the cardinality of  $\mathfrak{B}_2$ is 
$$|E^P\setminus K^P|=|E^P|-|K^P|.$$

Hence we have
\begin{align*}
\text{dim}(H^{*}(P_n)^{\mathfrak{G}})&=|\mathfrak{B}_1|+|\mathfrak{B}_2|-|\mathfrak{B}_1\cap\mathfrak{B}_2|\\
&=\frac{\text{dim}(H^{*}(P_n)^{\mathfrak{S}_{q} \times \mathfrak{S}_q})-|E^P|}{2}+|E^P|-|K^P|\\
&=\frac{\text{dim}(H^{*}(P_n)^{\mathfrak{S}_{q} \times \mathfrak{S}_q})+|E^P|}{2}-|K^P|.
\end{align*}

These prove the theorem. \qb

{\noindent\bf Theorem 4.3.}
{\it
$$|E^P|=\sum_{\lambda\in\Lambda(n-*)}
\sum_{(\lambda,d,k)\in E_\lambda}(\Pi_1\Pi_2\Pi_3)\:
$$
and
$$|K^P|=
\sum_{\lambda\in\Lambda(n-*)}
\sum_{(\lambda,d,k)\in K_\lambda}(\Pi_1\Pi_2\Pi_3),
$$
where
\begin{align*}
	\Pi_1&=\prod_{(\lambda_i,d_i)\in E_1(\lambda,d)}
	(\sum_{b=1}^{|(\lambda_i,d_i)|}|M( b)|\left(\begin{array}{c}|\Pi(\lambda_i,d_i)|\\b\end{array}\right)),\\
	\Pi_2&=\prod_{(\lambda_i,d_i)\in E_{21}(\lambda,d,k)}\left(\begin{array}{c}|\Pi(\lambda_i,d_i)|\\|(\lambda,d)_i|\end{array}\right)\prod_{(\lambda_i,d_i)\in E_{22}(\lambda,d,k)}
	\left(\begin{array}{c}\frac{|\Pi(\lambda_i,\frac{\lambda_i}{2})|-|\overline{\Pi(\lambda_i,\frac{\lambda_i}{2})|}}{2}\\|(\lambda,d)_i|\end{array}\right),\\
	\Pi_3&=\prod_{(\lambda_i,d_i)\in E_{3}(\lambda,d,k)}\left(\begin{array}{c}|\overline{\Pi(\lambda_i,\frac{\lambda_i}{2})}|\\|(\lambda,d)_i|\end{array}\right)
\end{align*}
and
 $\left( \begin{array}{c} 
	* \\ 
	* 
\end{array} \right)$ is the binomial coefficient.}

{\noindent\bf Proof.}
Let $$\lambda=(\underbrace{\lambda_1,\cdots,\lambda_1}_{|(\lambda_1)|},\underbrace{\lambda_2,\cdots,\lambda_2}_{|(\lambda_2)|},\cdots,\underbrace{\lambda_a,\cdots,\lambda_a}_{|(\lambda_a)|},1,\cdots,1)\in\Lambda(n-*).$$
We fix a $(\lambda,d,k) \in E_\lambda$. Recall that $s\otimes \zeta_{\lambda}\in E^P$ or $s\otimes \zeta_{\lambda}\in K^P$  is uniquely determined by $$\chi_{l_j+i}(\delta_s),$$ $1\leq i \leq |(\lambda_{j})|-k_j$ for each $\lambda_j$, $1\leq j\leq a$. Hence we only need to consider those $1\leq i\leq l_{a+1}$ such that $d_i\geq \frac{\lambda_j}{2}$, $1\leq j\leq a$.

If $(\lambda_i,d_i)\in E_1(\lambda,d)$, then for a fixed $(\lambda_i,d_i)$, $1\leq b \leq |(\lambda,d)_i|$ and  $(m_1,\cdots,m_b)\in M(b)$, we can fix some $s\otimes \zeta_{\lambda}\in E^P$ such that 
\begin{align*}
	\chi_{j}(\delta_s)&=\chi_{j+1}(\delta_s)=\cdots=\chi_{j+m_1-1}(\delta_s),\\
	\chi_{j+m_1}(\delta_s)&=\chi_{j+m_1+1}(\delta_s)=\cdots=\chi_{j+m_1+m_2-1}(\delta_s),\\
	&\cdots\\
	\chi_{j+m_1+\cdots+m_{b-1}}(\delta_s)&=\chi_{j+m_1+\cdots+m_{b-1}+1}(\delta_s)=\cdots=\chi_{j+|(\lambda,d)_i|-1}(\delta_s),
\end{align*}
where $j$ is fixed by $(\lambda_i,d_i)$ and $$\chi_{j}(\delta_s)<\chi_{j+m_1}(\delta_s)<\cdots<\chi_{j+m_1+\cdots+m_{b-1}}\in \Pi(\lambda_i,d_i).$$
It follows that the number of those elements stated above is $\left(\begin{array}{c}|\Pi(\lambda_i,d_i)|\\b\end{array}\right)$.

If $(\lambda_i,d_i)\in E_{21}(\lambda,d,k)$, then Theorem $2.3$ $(4)$ implies that for a fixed $(\lambda_i,d_i)$, we can fix some $s\otimes \zeta_{\lambda}\in E^P$ such that 
$$\chi_j(\delta_s)<\chi_{j+1}(\delta_s)<\cdots<\chi_{j+|(\lambda,d)_i|-1}(\delta_s) \in\Pi(\lambda_i,d_i),$$ where $j$ is fixed by $(\lambda_i,d_i)$. It follows that the number of those elements stated above is $\left(\begin{array}{c}|\Pi(\lambda_i,d_i)|\\|(\lambda,d)_i|\end{array}\right)$.

If $(\lambda_i,d_i)\in E_{22}(\lambda,d,k)$, we recall that  $d_i=\frac{\lambda_i}{2}$, then Theorem $2.3$ $(4)$ implies that for a fixed $(\lambda_i,d_i)$, we can fix some $s\otimes \zeta_{\lambda}\in E^P$ such that 
$$\chi_j(\delta_s)<\chi_{j+1}(\delta_s)<\cdots<\chi_{j+|(\lambda,d)_i|-1}(\delta_s)$$
and
$$\chi_{j+l}(\delta_s)<\overline{\chi_{j+l}}(\delta_s), 0\leq l\leq |(\lambda,d)_i|-1,$$
 where $j$ is fixed by $(\lambda_i,d_i)$. Note that for such $\chi_{j+l}(\delta_s)$, we have $$\chi_{j+l}(\delta_s)\in\left(\Pi(\lambda_{i},\frac{\lambda_i}{2})\setminus\overline{\Pi(\lambda_{i},\frac{\lambda_i}{2})}\right)/\langle\sigma\rangle.$$
 It follows that the number of those elements stated above is
 $$\left(\begin{array}{c}\frac{|\Pi(\lambda_i,\frac{\lambda_i}{2})|-|\overline{\Pi(\lambda_i,\frac{\lambda_i}{2})|}}{2}\\|(\lambda,d)_i|\end{array}\right).$$
 
 If $(\lambda_i,d_i)\in E_{3}(\lambda,d,k)$, then Theorem $2.3$ $(4)$ implies that for a fixed $(\lambda_i,d_i)$, we can fix some $s\otimes \zeta_{\lambda}\in E^P$ such that 
 $$\chi_j(\delta_s)<\chi_{j+1}(\delta_s)<\cdots<\chi_{j+|(\lambda,d)_i|-1}(\delta_s)$$
 and
 $$\chi_{j+l}(\delta_s)=\overline{\chi_{j+l}}(\delta_s), 0\leq l\leq |(\lambda,d)_i|-1,$$
 where $j$ is fixed by $(\lambda_i,d_i)$. It follows that the number of those elements stated above is
 $\left(\begin{array}{c}|\overline{\Pi(\lambda_i,\frac{\lambda_i}{2})}|\\|(\lambda,d)_i|\end{array}\right)$.

 Finally, since $(m_1,\cdots,m_b)$ runs over $M(b)$, $b$ runs over $1$ to $|(\lambda,d)_i|$, $i$ runs over $E_{1}(\lambda,d)\cup E_{21}(\lambda,d,k)\cup E_{22}(\lambda,d,k)\cup E_{3}(\lambda,d,k)$ and $(\lambda,d)$ runs over $E_\lambda$ or $K_\lambda$, $\lambda$ runs over $\Lambda(n-*)$, we get these formulas  consequently. \qb

Theorem $3.4$ in \cite{W24} directly provides the cardinality of $\Pi(\lambda_i,d_i)$.  Now we proceed to compute the cardinality of $\overline{\Pi(\lambda_i,\frac{\lambda_i}{2})}=\overline{\Pi(2d,d)}$ for some $d\in \mathbb{Z}_+$. By the direct computation, we have $\overline{\Pi(2,1)}=\{(1)\}$.

Now we suppose $2d\geq3$.
We recall that if $\chi_i(\delta_s)=\overline{\chi_i}(\delta_s)$, then  there exists some $1\leq i_0\leq m$, such that
$$\chi_i(\delta_s)=(\underbrace{0,\cdots,0}_{a_{m}-1},a_{i_0+1},\cdots,\underbrace{0,\cdots,0}_{a_{i_0+1}-1},a_{m},\underbrace{0,\cdots,0}_{a_{i_0}-1},a_{1},\cdots\underbrace{0,\cdots,0}_{a_1-1},a_{i_0}).$$

Note that as a sequence, for a fixed $1\leq i_0\leq m$, if there exists some $$c_d^{a_m+a_{m-1}+\cdots+a_{m-i+1}},$$ $1\leq i \leq m$, such that $$\begin{array}{rl}
	&c_d^{a_m+a_{m-1}+\cdots+a_{m-i+1}} \circ \chi_i(\delta_s)\\
	=&(\underbrace{0,\cdots,0}_{a_{m-i}-1},a_{i_0+i+1},\cdots,,\underbrace{0,\cdots,0}_{a_{i_0+1}-1},a_{m},\underbrace{0,\cdots,0}_{a_{i_0}-1},a_{1},\cdots,\underbrace{0,\cdots,0}_{a_{m-i+1}-1},a_{i_0+i})\\
	=&(\underbrace{0,\cdots,0}_{a_{m}-1},a_{i_0+1},\cdots,\underbrace{0,\cdots,0}_{a_{i_0+1}-1},a_{m},\underbrace{0,\cdots,0}_{a_{i_0}-1},a_{1},\cdots\underbrace{0,\cdots,0}_{a_1-1},a_{i_0})\\
	=&\chi_i(\delta_s),\end{array}$$
then for such $i_0$,  the sequences $(a_{i_0+1},\cdots,a_m,a_1,\cdots,a_{i_0})$ and $(a_m-1,\cdots,a_1-1)$ satisfy $$\begin{array}{rl}
&c_m^{i}\circ(a_{i_0+1},\cdots,a_m,a_1,\cdots,a_{i_0})\\
=&(a_{i_0+i+1},\cdots,a_m,a_1,\cdots,a_{i_0+i})\\
=&(a_{i_0+1},\cdots,a_m,a_1,\cdots,a_{i_0})
\end{array}$$
and
$$\begin{array}{rl}
	&c_m^{i}\circ(a_m-1,\cdots,a_1-1)\\
	=&(a_{m-i}-1,\cdots,a_1-1,a_m,\cdots,a_{m-i+1}-1)\\
	=&(a_m-1,\cdots,a_1-1)
\end{array},$$
where $c_m^{i+1}\in\langle c_m \rangle$ is an element in the cyclic group generated by the cycle $c_m$ of length $m$. Hence the order of  $\langle c_d \rangle_{\chi_i(\delta_s)}$ equals the order of  $\langle c_m \rangle_{(a_1,a_{2},\cdots,a_m)}$,
where we recall that $\langle c_{d}\rangle_{\chi_i(\delta_s)}$ is the isotropy group of the $\langle c_{d}\rangle$ action on  $\chi_{i}(\delta_{s})$ and $\langle c_m \rangle_{(a_1,a_{2},\cdots,a_m)}$ is the isotropy group of the $\langle c_{m}\rangle$ action on  $(a_1,a_{2},\cdots,a_m)$. 	
On the other hand, for a fixed $(a_1,a_{2},\cdots,a_m)$, for every $1\leq i_0\leq m$, there exists $$\chi_i(\delta_s)=(\underbrace{0,\cdots,0}_{a_{i_0}-1},a_{1},\cdots\underbrace{0,\cdots,0}_{a_1-1},a_{i_0},\underbrace{0,\cdots,0}_{a_{m}-1},a_{i_0+1},\cdots,\underbrace{0,\cdots,0}_{a_{i_0+1}-1},a_{m})$$
such that $\chi_i(\delta_s)=\overline{\chi_i}(\delta_s)$. Hence for a fixed $(a_1,a_{2},\cdots,a_m)$, there exists at most $m$ invariant cycles $\chi_i(\delta_s)$ such that  the order of  $\langle c_d \rangle_{\chi_i(\delta_s)}$ equals the order of  $\langle c_m \rangle_{(a_1,a_{2},\cdots,a_m)}$.

 We consider
\begin{align*}
	t&=(t_1,t_2,\cdots,t_p)\in(\mathbb{Z}_{>0})^p,\\
	a&=(a_{i_1}<a_{i_2}<\cdots<a_{i_p})\in(\mathbb{Z}_{
		>0})^p
\end{align*}
and the pair $T=(t,a)$.
It is easy to see that if the pair $T$ satisfies $$t_1+t_2+\cdots+t_p=m$$ and
$$t_1a_{i_1}+t_2a_{i_2}+\cdots+t_pa_{i_p}=d,$$ then there exist some $(a_1,a_{2},\cdots,a_m)$ such that
$$\{a_1,a_{2},\cdots,a_m\}=\{a_{i_1}<a_{i_2}<\cdots<a_{i_p}\},$$
where $a_{i_j}$ repeated $t_j$ times in $(a_1,a_{2},\cdots,a_m)$.
We write $T(p)$ the set of all $T$ which satisfy the above conditions.

For  $T=(t,a)\in T(p)$ and a divisor $k \geq 1$ of $gcd(t_1,t_2,\cdots,t_p)$, we write $$gcd(t_1,t_2,\cdots,t_p)=kp_1^{i_1}p_2^{i_2}\cdots p_s^{i_s},$$
where $p_j$, $1\leq j\leq s$, are distinct prime numbers.
Then we denote $$P(T,k)=\{p_1,p_2,\cdots,p_s\}.$$

We denote  $\left(\begin{array}{c}
	n\\
	k_1,
	k_2,
	\cdots,
	k_i
\end{array}\right)$ the multinomial coefficient, $k_1+k_2+\cdots+k_i=n$.

{\noindent\bf Definition 4.3.}(\cite{W24} Definition 3.3.)
{\it For $0 \leq j \leq s$, we denote by $P(T,k)_j$ the set consisting of all subsets of $P(T,k)$ with cardinality $j$.
	We agree $P(T,k)_0=\{\emptyset\}$.
	For a  divisor $k\geq 1$ of $gcd(t_1,t_2,\cdots,t_p)$ and $0\leq j \leq s$, we define the number $\pi_k(T,j)$ by
	\begin{align*}\pi_k(T,j)&=\sum_{\{p_{x_1},p_{x_2},\cdots, p_{x_j}\}\in P(T,k)_{j}}\left(\begin{array}{c}
			\frac{m}{kp_{x_1}p_{x_2}\cdots p_{x_j}}\\
			\frac{t_1}{kp_{x_1}p_{x_2}\cdots p_{x_j}},
			\frac{t_2}{kp_{x_1}p_{x_2}\cdots p_{x_j}},
			\cdots,
			\frac{t_p}{kp_{x_1}p_{x_2}\cdots p_{x_j}}
		\end{array}\right),
	\end{align*}
	while 
	\begin{align*}
		\pi_k(T,0)&=\left(\begin{array}{c}
			\frac{m}{k}\\
			\frac{t_1}{k},
			\frac{t_2}{k},
			\cdots,
			\frac{t_p}{k}
		\end{array}\right).
	\end{align*}
	We set $\pi_k(T,j)=0$ for every $0 \leq j \leq s$ if $k$ is not a divisor of $gcd(t_1,t_2,\cdots,t_p)$.
}

{\noindent\bf Theorem 4.4.}
{\it
	\begin{enumerate}
		\item[\bf 4.4.1.] If  $2d\equiv 0$ mod $4$, then $$|\overline{\Pi(2d,d)}|=\sum_{m=1}^{d}\sum_{p=1}^{m}\sum_{T\in T(p)}(\sum_{j=0}^{s}(-1)^j\pi_1(T,j)).$$
		
		\item[\bf 4.4.2.] If  $2d\equiv2$ mod $4$, then
		$$
		|\overline{\Pi(2d,d)}|=\sum_{m=1}^{d}\sum_{p=1}^{m}\sum_{T\in T(p)}\left((\sum_{j=0}^{s}(-1)^j\pi_1(T,j))
		+(\sum_{j=0}^{s}(-1)^j\pi_2(T,j))\right).
		$$
	\end{enumerate}
}

{\noindent\bf Proof.}  For a fixed $1\leq m\leq d$ and a fixed $T=(t,a)\in T(p)$, let $(a_1,a_{2},\cdots,a_m)\in(\mathbb{Z}_{>0})^{m}$ such that $$\{a_1,a_{2},\cdots,a_m\}=\{a_{i_1}<a_{i_2}<\cdots<a_{i_p}\},$$
$$t_1+t_2+\cdots+t_p=m$$ and $$t_1a_{i_1}+t_2a_{i_2}+\cdots+t_pa_{i_p}=d.$$ Then the order of 
$\langle c_m\rangle_{(a_1,a_{2},\cdots,a_m)}$  is a divisor of $gcd(t_1,t_2,\cdots,t_p)$.

We denote by $\mathfrak{O}_0$ the set of those elements  $$(a_1,a_{2},\cdots,a_m)\in(\mathbb{Z}_{>0})^{m}$$ such that $\langle c_m\rangle_{(a_1,a_{2},\cdots,a_m)}$ is of order $k$,
$\mathfrak{O}_i$ the set of those elements  $$(a_1,a_{2},\cdots,a_m)\in(\mathbb{Z}_{>0})^{m}$$ such that the order of $\langle c_m\rangle_{(a_1,a_{2},\cdots,a_m)}$  equals a positive multiple of $kp_i$.
Then we have $$\bigcup_{i=0}^{s}\mathfrak{O}_i$$ is the set of those elements  $$(a_1,a_{2},\cdots,a_m)\in(\mathbb{Z}_{>0})^{m}$$ such that the order of $\langle c_m\rangle_{(a_1,a_{2},\cdots,a_m)}$   equals a positive multiple of $k$.

Assume that $$(a_1,a_{2},\cdots,a_m)\in (\mathbb{Z}_{>0})^{\frac{m}{k}}$$ satisfies
$$\{a_1,a_2,\cdots,a_{\frac{m}{k}}\}=\{a_{i_1},a_{i_2},\cdots,a_{i_p}\},$$ $\#(a_{i_j})=\frac{t_j}{k}$, $1\leq j \leq p$. If we extend it to $$(a_1,a_2,\cdots,a_{\frac{m}{k}},a_1,a_2,\cdots,a_{\frac{m}{k}},\cdots,a_1,a_2,\cdots,a_{\frac{m}{k}})=(a_1,a_2,\cdots,a_m),$$ then it is an element in $(\mathbb{Z}_{>0})^{m}$ such that
the order of $\langle c_m\rangle_{(a_1,a_2,\cdots,a_m)}$ equals a positive multiple of $k$, which depends on the choice of $$(a_1,a_2,\cdots,a_{\frac{m}{k}}).$$ Since $\pi_k(T,0)$ is the number of all such elements, we have $$\pi_k(T,0)=|\bigcup_{i=0}^{s}\mathfrak{O}_i|.$$

Similarly, for a fixed $\{p_{x_1},p_{x_2},\cdots, p_{x_j}\}\in P(T,k)_j$, assume that $$(a_1,a_2,\cdots,a_{\frac{m}{kp_{x_1}p_{x_2}\cdots p_{x_j}}})\in(\mathbb{Z}_{>0})^{\frac{m}{kp_{x_1}p_{x_2}\cdots p_{x_j}}}$$ satisfies
$$\{a_1,a_2,\cdots,a_{\frac{m}{kp_{x_1}p_{x_2}\cdots p_{x_j}}}\}=\{a_{i_1},a_{i_2},\cdots,a_{i_p}\},$$ $\#(a_{i_j})=\frac{t_j}{kp_{x_1}p_{x_2}\cdots p_{x_j}}$, $1\leq j \leq p$.  If we extend it to
\begin{align*}&(a_1,a_2,\cdots,a_{\frac{m}{kp_{x_1}p_{x_2}\cdots p_{x_j}}},\cdots,a_1,a_2,\cdots,a_{\frac{m}{kp_{x_1}p_{x_2}\cdots p_{x_j}}})\\
	=&(a_1,a_2,\cdots,a_m),\end{align*} then it is an element in $(\mathbb{Z}_{>0})^{m}$ such that
the order of $\langle c_m\rangle_{(a_1,a_2,\cdots,a_m)}$ equals a positive multiple of $kp_{x_1}p_{x_2}\cdots p_{x_j}$, which depends on the choice of $$(a_1,a_2,\cdots,a_{\frac{m}{kp_{x_1}p_{x_2}\cdots p_{x_j}}}).$$ Since the summand  of $\pi_k(T,j)$ fixed by $\{p_{x_1},p_{x_2},\cdots, p_{x_j}\}\in P(T,k)_j$ is the number of all such elements,  we have $$\pi_k(T,j)=\sum_{\{p_{x_1},p_{x_2},\cdots ,p_{x_j}\}\in P(T,k)_j}|\bigcap_{i=1}^{j}\mathfrak{O}_{x_i}|.$$

Since $\mathfrak{O}_0 \cap \mathfrak{O}_i=\emptyset$, $1\leq i\leq s$,
by applying the inclusion-exclusion principle, we have
$$\pi_k(T,0)=|\mathfrak{O}_0|-\sum_{j=1}^{s}(-1)^j\pi_k(T,j).$$
Thus $$|\mathfrak{O}_0|=\sum_{j=0}^{s}(-1)^j\pi_k(T,j)$$ is the number of all
$$(a_1,a_2,\cdots,a_m)\in(\mathbb{Z}_{>0})^{m}$$ such that  $\langle c_m\rangle_{(a_1,a_2,\cdots,a_m)}$ is of order $k$ and
$\frac{k}{m}|\mathfrak{O}_0|$ is the number of all $$(a_1,a_2,\cdots,a_m)\in(\mathbb{Z}_{>0})^{m}/\langle c_m\rangle$$ such that $\langle c_m\rangle_{(a_1,a_2,\cdots,a_m)}$ is of order $k$.

On the other hand, for a fixed $$(a_1,a_2,\cdots,a_{\frac{m}{k}},\cdots,a_1,a_2,\cdots,a_{\frac{m}{k}})=(a_1,a_2,\cdots,a_m)\in (\mathbb{Z}_{>0})^{m}/\langle c_m\rangle$$ and a fixed $1\leq i_0 \leq m$,
\begin{align*}
	&\chi_i(\delta_s)\\
	=&(\underbrace{0,\cdots,0}_{a_{i_0}-1},a_{1},\underbrace{0,\cdots,0}_{a_{i_0-1}-1},a_{2},\cdots,\underbrace{0,\cdots,0}_{a_{i_0+1}-1},a_{\frac{m}{k}},\cdots,\underbrace{0,\cdots,0}_{a_{i_0}-1},a_{1},\cdots,\underbrace{0,\cdots,0}_{a_{i_0+1}-1},a_{\frac{m}{k}})
\end{align*}
is a invariant cycle such that $\langle c_d \rangle_{\chi_i(\delta_s)}$ is of order $k$. If $1\leq i_0\neq i_0'\leq \frac{m}{k}$, then
\begin{align*}
	&(\underbrace{0,\cdots,0}_{a_{i_0}-1},a_{1},\underbrace{0,\cdots,0}_{a_{i_0-1}-1},a_{2},\cdots,\underbrace{0,\cdots,0}_{a_{i_0+1}-1},a_{\frac{m}{k}},\cdots,\underbrace{0,\cdots,0}_{a_{i_0}-1},a_{1},\cdots,\underbrace{0,\cdots,0}_{a_{i_0+1}-1},a_{\frac{m}{k}})\\
	\neq&(\underbrace{0,\cdots,0}_{a_{i_0'}-1},a_{1},\underbrace{0,\cdots,0}_{a_{i_0'-1}-1},a_{2},\cdots,\underbrace{0,\cdots,0}_{a_{i_0'+1}-1},a_{\frac{m}{k}},\cdots,\underbrace{0,\cdots,0}_{a_{i_0'}-1},a_{1},\cdots,\underbrace{0,\cdots,0}_{a_{i_0'+1}-1},a_{\frac{m}{k}}).
\end{align*}
While if $i_0\equiv i_0'$ mod $\frac{m}{k}$, then 
\begin{align*}
	&(\underbrace{0,\cdots,0}_{a_{i_0}-1},a_{1},\underbrace{0,\cdots,0}_{a_{i_0-1}-1},a_{2},\cdots,\underbrace{0,\cdots,0}_{a_{i_0+1}-1},a_{\frac{m}{k}},\cdots,\underbrace{0,\cdots,0}_{a_{i_0}-1},a_{1},\cdots,\underbrace{0,\cdots,0}_{a_{i_0+1}-1},a_{\frac{m}{k}})\\
	=&(\underbrace{0,\cdots,0}_{a_{i_0'}-1},a_{1},\underbrace{0,\cdots,0}_{a_{i_0'-1}-1},a_{2},\cdots,\underbrace{0,\cdots,0}_{a_{i_0'+1}-1},a_{\frac{m}{k}},\cdots,\underbrace{0,\cdots,0}_{a_{i_0'}-1},a_{1},\cdots,\underbrace{0,\cdots,0}_{a_{i_0'+1}-1},a_{\frac{m}{k}}).
\end{align*}
Hence the parameter
$i_0$ yields exactly $\frac{m}{k}$ distinct invariant cycles,  the number of all $\chi_i(\delta_s)$ such that $\langle c_d \rangle_{\chi_i(\delta_s)}$ is of order $k$ is $|\mathfrak{O}_0|$.

By applying Theorem $2.3$ $(2)$ and $(3)$,  the theorem follows. \qb

{\bf Remark:} In fact, Theorem $4.4$ is a corollary of Theorem $3.4$ in \cite{W24}. The difference is that we consider \begin{align*}
	t&=(t_1,t_2,\cdots,t_p)\in(\mathbb{Z}_{>0})^p,\\
	a&=(a_{i_1}<a_{i_2}<\cdots<a_{i_p})\in(\mathbb{Z}_{
		>0})^p
\end{align*} for some sequences $$(a_{1},\cdots,a_{m})\in(\mathbb{Z}_{>0})^{m}$$ such that $$t_1+t_2+\cdots+t_p=m$$ and
$$t_1a_{i_1}+t_2a_{i_2}+\cdots+t_pa_{i_p}=d$$ instead of
consider
\begin{align*}
	t&=(t_1,t_2,\cdots,t_p)\in(\mathbb{Z}_{>0})^p,\\
	a&=(a_{i_1}<a_{i_2}<\cdots<a_{i_p})\in(\mathbb{Z}_{\geq0})^p
\end{align*} for some invariant cycles $$\chi_{i}(\delta_s)=(a_{1},\cdots,a_{d})\in(\mathbb{Z}_{\geq0})^{d}$$
such that
$$t_1+t_2+\cdots+t_p=d$$ and
$$t_1a_{i_1}+t_2a_{i_2}+\cdots+t_pa_{i_p}=\lambda-d.$$

\end{document}